\documentclass[reqno]{amsart}
 \usepackage{amsbsy,amssymb,amscd,amsfonts,latexsym,amstext,delarray,
 amsmath,graphicx, color}
 \usepackage[dvips]{epsfig}
\input xypic
\usepackage{hyperref}

\definecolor{Green}{rgb}{0.0,0.40,0.0}

\newtheorem{thm}{Theorem}[section]
\newtheorem{prop}[thm]{Proposition}
\newtheorem{cor}[thm]{Corollary}
\newtheorem{lem}[thm]{Lemma}

\newtheorem{defn}[thm]{Definition}
\newtheorem{rem}[thm]{Remark}
\newtheorem{example}[thm]{Example}

\def\qqq{\,,\quad~\forall}

\def\Aut{{\rm Aut}}

\def\Spec{{\rm Spec}}

\def\B{{\mathbb B}}
\def\C{{\mathbb C}}
\def\F{{\mathbb F}}

\def\N{{\mathbb N}}

\def\Q{{\mathbb Q}}
\def\R{{\mathbb R}}
\def\Z{{\mathbb Z}}

\def\W{{\mathbb W}}

\def\cD{{\mathcal D}}

\def\cN{{\mathcal N}}

\def\cT{{\mathcal T}}
\def\cU{{\mathcal U}}

\newcommand{\ie}{{\it i.e.\/}\ }
\newcommand{\eg}{{\it e.g.\/}\ }
\newcommand{\cf}{{\it cf.\/}\ }

\def\rmax{\R^{\rm{max}}_+}

\def\te{Teichm\"uller }
\def\mc{multiplicatively cancellative }
\def\run{\R^{\rm un}}

\newcommand{\nil}[1]{}

\title
{The Witt construction in characteristic one  and Quantization}

\author[Connes]{Alain Connes}
\address{A.~Connes: Coll\`ege de France \\
3, rue d'Ulm \\ Paris, F-75005 France
\\ I.H.E.S. and Vanderbilt
University} \email{alain\@@connes.org}


\date{}
\begin{document}
\maketitle

\begin{center}
{\em  Dedicated to Henri Moscovici, with admiration and friendship}
\end{center}

\begin{abstract} We develop the analogue of the Witt construction in characteristic one. We construct a functor from pairs $(R,\rho)$ of a perfect semi-ring $R$ of characteristic one and an element $\rho>1$ of $R$ to real Banach algebras. We find that the entropy function occurs uniquely as the analogue of the \te polynomials in characteristic one. We then apply the construction to the semi-field $\rmax$ which plays a central role in idempotent analysis and tropical geometry. Our construction gives the inverse process of the ``dequantization" and provides a first hint towards an extension $\run$ of the field of real numbers relevant both in number theory and quantum physics.
\end{abstract}

\section{Foreword}
The celebration of Henri's 65th birthday gives me a long awaited occasion to
express my deep gratitude to him, for his indefectible friendship in our long journey through mathematics, since the first time we met in Princeton in the fall of 1978. Besides the great enlightening moments,
 those that I cherish most are the times when we both knew we were close to some ``real stuff" but also knew that we could get there only at the price of time consuming efforts which we shared so happily over the years.
\tableofcontents

\section{Introduction}

The goal of this paper is to develop an analogue of the Witt construction in the case of characteristic one, which was initiated in \cite{jamifine}. Our starting point is a formula which goes back to \te and which gives an explicit expression for the sum of the multiplicative lifts in the context of strict $p$-rings. A ring $R$ is a strict $p$-ring when $R$ is complete
and Hausdorff with respect to the $p$-adic topology, $p$ is not a zero-divisor in $R$, and the residue ring
$K = R/pR$ is perfect. The ring $R$ is uniquely determined by $K$ up to canonical isomorphism and there exists a unique multiplicative section $\tau:K\to R$ of the residue morphism $\epsilon: R \to K=R/pR$
\begin{equation}\label{3m}
\tau: K\to R,\quad \epsilon\circ\tau = id ,\ \tau(xy)=\tau(x)\tau(y)\qqq x,y\in K.
\end{equation}
 Every element $x$ of $R$ can be written uniquely in the form
\begin{equation}
x = \sum\tau(x_n)p^n\qquad x_n\in K
\end{equation}
which gives a canonical bijection $\tilde\tau: K[[T]]\to R$ such that
\begin{equation}\label{bijcan}
   \tilde\tau(\sum x_nT^n)=\sum\tau(x_n)p^n.
\end{equation}
The formula which goes back to \te \cite{Teich} allows one to express the sum of two (or more) multiplicative lifts in the form
\begin{equation}\label{sumtau}
 \tau(x)+ \tau(y) = \tilde\tau\left(\sum_{\alpha\in I_p}w_p(\alpha,T)\,x^\alpha y^{1-\alpha}\right)
\end{equation}
In this equation the sum inside the parenthesis in the right hand side takes place in $K[[T]]$, the variable $\alpha$ ranges in
\begin{equation}\label{20m}
I_p =\{ \alpha\in\Q\cap[0,1],\quad\exists n,~p^n\alpha\in\Z\}
\end{equation}
so that, since $K$ is perfect, the terms $x^\alpha y^{1-\alpha}$ make sense. Finally the terms
\begin{equation}\label{wp}
    w_p(\alpha,T)\in \F_p[[T]]\qqq \alpha \in I_p
\end{equation}
only depend on the prime $p$ and tend to zero at infinity in $\F_p[[T]]$, for the discrete topology in $I_p$, so that the sum \eqref{sumtau} is convergent. The formula \eqref{sumtau} easily extends to express the sum of $n$ multiplicative lifts as
\begin{equation}\label{sumtaun}
\sum \tau(x_j)  = \tilde\tau\left(\sum_{\alpha_j\in I_p,\,\sum\alpha_j=1}w_p(\alpha_1,\ldots,\alpha_n,T)\,\prod x_j^{\alpha_j}\right)
\end{equation}
As is well known the algebraic structure of $R$ was functorially reconstructed from that of $K$ by Witt who showed that the algebraic rules in $R$ are polynomial in terms of the components $X_n$
\begin{equation}\label{wittcomp}
 x = \sum\tau(X_n^{p^{-n}})p^n\,, \ \   X_n=x_n^{p^{n}}
\end{equation}
which makes sense since $K$ is perfect.
One can in fact also reconstruct the full algebraic structure of $R$ as a deformation of $K$ depending upon the parameter $T$ using the above formula \eqref{sumtaun} but the corresponding algebraic relations are not so simple to handle, mostly because the map $x\to x^\alpha$ from $K$ to $K$ is an automorphism
only when $\alpha$ is a  power of $p$. As we shall now explain, this difficulty disappears in the limit case of characteristic one. The framework we use for characteristic one is that of semi-rings
(\cf \cite{Golan}). For \mc semi-rings $S$ of characteristic one the maps $x\to x^n$ are injective endomorphisms for any positive integer $n$ (\cf  \cite{Golan}) and it is natural to say that $S$ is perfect when these maps are surjective. It then follows that the fractional powers $x\to x^\alpha$ make sense for any $\alpha\in \Q_+^*$ and define automorphisms $\theta_\alpha\in \Aut(S)$. We shall show in \S \ref{wittone} how to solve the functional equation on the coefficients $w(\alpha)$ defined for $\alpha\in \bar I=\Q\cap [0,1]$ which ensures  that the following analogue of \eqref{sumtau} defines a deformation of $S$ into a semi-ring of characteristic zero
\begin{equation}\label{deform}
    x+_w y=\sum_{\alpha\in \bar I}w(\alpha)\, x^\alpha y^{1-\alpha}
\end{equation}
We let $w(0)=w(1)=1$ and $I=\Q\cap (0,1)$ so that $\bar I=I\cup\{0,1\}$.  The commutativity means that
\begin{equation}\label{comm}
w(1-\alpha) = w(\alpha)\,,
\end{equation}
and the associativity means that the equality
\begin{equation}\label{dec1}
    w^{(3)}(\alpha_1,\alpha_2,\alpha_3)=w(\alpha_1)w(\frac{\alpha_2}{1-\alpha_1})^{1-\alpha_1} \qqq \alpha_j \in I \mid \sum\alpha_j=1
\end{equation}
defines a symmetric function on the simplex $\Sigma_3$ where, more generally,
\begin{equation}\label{simplex}
    \Sigma_n=\{(\alpha_1,\ldots,\alpha_n)\in I^{n}\mid \sum\alpha_j = 1\}
\end{equation}
Note that \eqref{comm} means that the function $w^{(2)}(\alpha_1,\alpha_2)=w(\alpha_1)$ is symmetric on $\Sigma_2$. The equations of symmetry of $w^{(n)}$, $n=2,3$ only use the multiplicative structure of $S$ and thus continue to make sense for any map $w:I\to G$ where $G$ is a uniquely divisible abelian group (denoted multiplicatively). We take $G=S^\times$   the multiplicative  group of invertible elements of $S$. We show in Theorem \ref{checkchar} that all solutions of these  equations (the symmetry of $w^{(n)}$, $n=2,3$) are of the form
\begin{equation}\label{allsol}
w(\alpha) = \chi(\alpha)^\alpha\chi(1-\alpha)^{1-\alpha}\qqq \alpha\in I
\end{equation}
where $\chi$ is a homomorphism $\Q_+^\times \to G$. Thus in this generality one can give arbitrarily the value of $\chi(p)\in G$ for all primes $p$. But the group of invertible elements of $S$ admits an additional structure: the  partial order coming from the additive structure of $S$. We show in Theorem \ref{entropythm} that, provided the $\theta_s$ extend by continuity to $s\in \R_+^*$,  all solutions
which fulfill the inequality
\begin{equation}\label{pos}
w(\alpha)\geq 1,\qquad\forall\alpha\in\Q\cap[0,1]\,,
\end{equation}
are of the form
\begin{equation}\label{entro}
w(\alpha)=\rho^{ S(\alpha)},~S(\alpha) = -\alpha\log(\alpha)-(1-\alpha)\log(1-\alpha).
\end{equation}
where  $\rho\in S$, $\rho\ge 1$. One recognizes the {\em entropy} function
\begin{equation}\label{entrop}
   S(\alpha) = -\alpha\log(\alpha)-(1-\alpha)\log(1-\alpha)
\end{equation}
which is familiar in thermodynamics, information theory and ergodic theory.
In \S \ref{wittonesect} we construct a functor $W$ from pairs $(R,\rho)$ of a \mc perfect semi-ring $R$ of characteristic one and an invertible element $\rho> 1$ in $R$ to algebras over $\R$. The construction of the algebra $W(R,\rho)$ involves several operations
 \begin{itemize}
   \item  A completion with respect to a $\rho$-adic distance canonically associated to $\rho$.
   \item A deformation of the addition involving the $w(\alpha)=\rho^{ S(\alpha)}$.
   \item A symmetrization to obtain a ring from a semi-ring.
 \end{itemize}
We also show in \S \ref{Banach} that the algebra $W(R,\rho)$ over $\R$ naturally yields a Banach algebra $\overline W(R,\rho)$ obtained by completion and still depending functorially on $(R,\rho)$.
The Gelfand spectrum of the complexified algebra $\overline W(R,\rho)_\C$ is a non-empty compact space canonically associated to $(R,\rho)$.

In the last section \ref{sectrun}, we return to a more algebraic set-up and to the analogy with the Witt construction in characteristic $p$. The need for the construction of an extension
$\run$ of $\R$ playing a role similar to the maximal unramified extension of the $p$-adic fields appears both in number theory and in quantum physics. In number theory we refer to the introduction of \cite{Weil} for the need of an interpretation of the connected component of identity in the id\`ele class group of the global  field of rational numbers as a Galois group involving a suitable refinement of the maximal abelian extension of $\Q$ (\cf also the last section of \cite{wagner}). This global question admits a local analogue whose solution requires constructing $\run$.   We show in \S \ref{sectrun} that the analogue of the Witt construction in characteristic one gives a first hint of what $\run$ could look like. One obtains the (completion of the) maximal unramified extension of $p$-adic integers by applying the Witt construction $\W_p$ (at the prime $p$) to the extension of the finite field $\F_p$ given by an algebraic closure $\bar\F_p$. The analogue, in characteristic one, of the latter extension is the extension
\begin{equation}\label{Br}
\B\subset \rmax
\end{equation}
of the only finite semi-field $\B$ which is not a field, by the semi-field $\rmax$ which plays a central role in idempotent analysis (\cf \cite{Maslov}, \cite{Litvinov}) and tropical geometry (\cite{Gathmann}, \cite{Mikhalkin}). The semi-field $\rmax$ is defined as the locally compact space $\R_+=[0,\infty)$ endowed with the ordinary product and a new idempotent addition $x +' y=\max\{ x ,y \}$. It admits a one parameter group of automorphisms $\theta_\lambda\in \Aut(\rmax)$, $\theta_\lambda(x)=x^\lambda$ for all $x\in\rmax$, which is the analogue of the arithmetic Frobenius. Since the above construction of $W(R,\rho)$ depends upon the choice of the element $\rho\in R$, $\rho>1$, the $W(\theta_\lambda)$ give canonical isomorphisms
\begin{equation}\label{shift}
W(\theta_\lambda): W(R,\rho)\to W(R,\theta_\lambda(\rho))
\end{equation}
To eliminate the dependence upon $\rho$ in the case of $\rmax$, we consider all $\rho$'s simultaneously, \ie we let $\rho=e^T$ where $T>0$ is a parameter. One then has
\begin{equation}\label{Tdepend}
w(\alpha,T) = \alpha^{-T\alpha}(1-\alpha)^{-T(1-\alpha)}
\end{equation}
which depends on the parameter $T$ a bit like the $w_p(\alpha,T)$ of \eqref{wp}. When one computes using the formula \eqref{deform}, \ie
\begin{equation}\label{deformT}
x+_w y=\sum_{\alpha\in I}w(\alpha,T)\, x^\alpha y^{1-\alpha}
\end{equation}
 the sum of $n$ terms $x_j$ independent of $T$  one obtains
\begin{equation}\label{sumofn}
    x_1+_w\cdots +_w x_n=\left(\sum x_j^{1/T}\right)^T
\end{equation}
In particular if all $x_j=1$, this gives the function $T\mapsto n^{T}$.  In fact more generally the functions of the form $T\mapsto x^{T}$ form a sub semi-ring isomorphic to $\R_+$ with its ordinary operations (of addition and multiplication) and are the fixed points of the natural extension of the $\theta_\lambda$ as automorphisms given, using \eqref{shift}, by
\begin{equation}\label{autos}
\alpha_\lambda(f)(T)=\theta_\lambda(f(T/\lambda))=f(T/\lambda)^\lambda\,.
\end{equation}
In general for elements given by functions $f(T)$ the addition coming from \eqref{deformT} is given by
\begin{equation}\label{addintro}
    (f+_wg)(T)=\left(f(T)^{1/T}+g(T)^{1/T}\right)^T\qqq T.
\end{equation}
This shows that for each $T$ there is a uniquely associated character $\chi_T$ with values in real numbers with their usual operations,  which is given by
\begin{equation}\label{charach1}
    \chi_T(f)=f(T)^{1/T}
\end{equation}
and one can use the characters $\chi_T$ to represent the elements of the extension $\run$ as functions of $T$ with the ordinary  operations of pointwise sum and product. In this representation the functions $\tau(x)$ which were independent of $T$ now are represented by
$$
T\mapsto \chi_T(\tau(x))=x^{1/T}
$$
Such functions $\tau(x)$ are the analogues of the \te lifts and they generate  the field formed of fractions of the form
\begin{equation}\label{simplest1}
\chi_T(f)=\left(\sum \lambda_j e^{-s_j/T}\right)\text{\huge/}\left(\sum \mu_j e^{-t_j/T}\right)
\end{equation}
where the coefficients $\lambda,\mu$ are rational numbers and the exponents $s,t$ are real numbers.
 Such fractions, with the coefficients $\lambda,\mu$ real, give a first hint towards $\run$ but one  should relax the requirement that the sums involved are finite. We briefly discuss in \S \ref{charep} the role of divergent series (\cite{Ramis}) in the representation of elements of $\run$. The key examples to be covered come  from quantum physics and are given by   functional integrals which are typically of the form (\cf \S \ref{deforsect})
\begin{equation}\label{ZJE}
 Z(J) =  \left(\int \exp (-\frac{S(\phi)- \langle J, \phi
 \rangle}{\hbar})\cD[\phi]\right)\text{\huge/}\left(\int \exp(-\frac{S(\phi)}{\hbar})\,
 \cD[\phi]\right).
\end{equation}
This suggests that, in quantum physics, the parameter $T$ should be related to $\hbar$
	 and the above expression \eqref{ZJE} should define for each value of the source parameter $J$ an element $f\in \run$ such that $Z(J) =\chi_T(f)$ for  $\hbar=T$. Moreover the elements of $\run$ obtained in this manner have asymptotic expansions of the form
 \begin{equation}\label{expexp}
 f(T)=\chi_T(f)^T\sim \sum a_nT^n\,,
 \end{equation}
where one sets $T=\hbar$. This property should be satisfied by all elements of $\run$ and is compatible with the fact that  the ``numbers" which are obtained after quantization can often only be described as sums of formal perturbative series in powers of $\hbar$. Basing the definition of $\run$ on the asymptotic expansion \eqref{expexp} is however loosing too much information as shown in the simple case \eqref{simplest1} where the expansion only depends upon the terms with lowest values of $s_j$ and $t_j$. We briefly relate this issue with the theory of Borel summation of divergent series in \S \ref{charep} and then show, at the end of \S \ref{sectrun}, that the above construction is intimately related to idempotent analysis (\cf \cite{Maslov}, \cite{Litvinov}). In fact our analogue in characteristic one of the Witt construction
provides the inverse  process of the ``dequantization" of idempotent analysis, a fact which justifies the word ``quantization" appearing in the title of this paper.

In conclusion we conjecture that the extension $\run$ of $\R$ is the natural home for the ``values" of the many $\hbar$-dependent physical quantities arising in quantum field theory. This fits  with the previous understanding of renormalization from the Riemann-Hilbert correspondence (\cf \cite{CK}, \cite{cknew}, \cite{cmln}, \cite{CMbook}).

\section{Sum of Teichm\"uller representatives}\label{wittcase}

Our goal in this section is to recall a formula which goes back to Teichm\"uller \cite{Teich}, for the sum of \te representatives in the context of strict $p$-rings.
We begin by recalling briefly the simplest instance of the polynomials with integral coefficients which express the addition in the Witt theory (We refer to \cite{Serre} and \cite{Rabi} for a quick introduction to that theory).
One defines the polynomials with integral coefficients $S_n(x,y)\in \Z[x,y]$ by the equality
\begin{equation}\label{1m}
(1- tx)(1-ty)=\prod_{n=1}^\infty(1-S_n(x,y)t^n)
\end{equation}
The first few are of the form
$$
\begin{array}{cc}
 S_1= & x+y \\
 S_2= & -x y \\
 S_3= & -x y (x+y) \\
 S_4= & -x y (x+y)^2 \\
 S_5= & -x y \left(x^3+2 x^2 y+2 x y^2+y^3\right) \\
 S_6= & -x y (x+y)^2 \left(x^2+x y+y^2\right) \\
 S_7= & -x y (x+y) \left(x^2+x y+y^2\right)^2 \\
 S_8= & -x y (x+y)^2 \left(x^4+2 x^3 y+4 x^2 y^2+2 x y^3+y^4\right) \\
 S_9= & -x y (x+y) \left(x^2+x y+y^2\right)^3 \\
 S_{10}= & -x y (x+y)^2 \left(x^6+3 x^5 y+7 x^4 y^2+8 x^3 y^3+7 x^2 y^4+3 x y^5+y^6\right)
\end{array}
$$
 One has for all n
\begin{equation}\label{9m}
x^n +y^n  = \sum_{d|n}d\,S_d(x,y)^{\frac nd}
\end{equation}
as follows by equating the coefficients of $t^n$ in $-t\frac{d}{dt}\log$ of both sides of \eqref{1m}.

\subsection{Strict $p$-rings}

We now fix a prime $p$ and recall well known notions.

\begin{defn}\label{pring}
 A ring $R$ is called a strict $p$-ring provided that $R$ is complete
and Hausdorff with respect to the $p$-adic topology, $p$ is not a zero-divisor in $R$, and the residue ring
$S = R/pR$ is perfect.
\end{defn}

Let $S$ be a perfect ring of characteristic $p$.

1. There is a strict $p$-ring $R$ with residue ring $S$, which is unique up to canonical isomorphism.

2. There exists a unique multiplicative section (called the \te section),
\begin{equation}\label{3mbis}
\tau: S\to R\quad \epsilon\circ\tau = id\qquad \epsilon: R \to S=R/pR
\end{equation}
3. Every element $x$ of $R$ can be written uniquely in the form
\begin{equation}
x = \sum\tau(x_n)p^n\qquad x_n\in K
\end{equation}
4. The construction of $R$ and $\tau$ is functorial in $S$: if $f:S\to S'$  is a homomorphism of
perfect rings of characteristic $p$,
then there is a unique homomorphism $F:R\to R'$ given by
\begin{equation}\label{functoria}
    F(\sum\tau(x_n)p^n)=\sum\tau(f(x_n))p^n
\end{equation}

\subsection{Universal coefficients $w_p(\alpha)$}

We consider a strict $p$-ring $R$ with a perfect residue ring $S = R/pR$. We let $\tau: S \to R$ be the multiplicative  section. Since $S$ is of characteristic $p$ and is perfect, we can define
\begin{equation}\label{12m}
\theta_\alpha(x) = x^\alpha\qquad\forall x\in S, ~\alpha = a/p^n\quad a\in\N, ~n\in\N
\end{equation}
Note that this is an automorphism of $S$ when $\alpha$ is a power of $p$ but not in general.
\begin{lem} With $S_k\in\Z[x,y]$ defined above, one has for any pair $(R,S)$ and any $x,y\in S$
\begin{equation}\label{13m}
\tau(x)+\tau(y) = \sum_0^\infty\tau(S_{p^n}(x^{1/p^n},y^{1/p^n}))p^n
\end{equation}
\end{lem}
\proof This follows from \cite{Rabi} Theorem 1.5 and Lemmas 3.2 and 3.4 which show that the polynomials $S_{p^n}(x,y)$ used above agree with the $S_n(x,y,0,\ldots,0)$ of \cite{Rabi} Theorem 1.5.
\endproof

We now expand the polynomial $S_{p^n}(x,y)$ in the form
\begin{equation}\label{19m}
S_{p^n}(x,y) = \sum_{k=0}^{p^n}a(n,k)x^ky^{p^n-k},\qquad a(n,k)\in\Z
\end{equation}
and we define a map $w_p$ from the set $I_p$ of \eqref{20m},
\begin{equation}\label{20mbis}
I_p = \{\alpha\in\Q\cap[0,1],\quad\exists n,~p^n\alpha\in\Z\}
\end{equation}
to the local ring $\F_p[[T]]$ of formal power series in $T$, by
\begin{equation}\label{21m}
w_p: I_p\to \F_p[[T]],\quad w_p(\alpha)=\sum_{k/p^n=\alpha}\bar a(n,k)T^n\quad \bar a(n,k) \equiv a(n,k)~\text{mod}~p
\end{equation}
\begin{lem} $a)$ One has $w_p(\alpha)\in\F_p[[T]]$ for all $\alpha\in I_p$.\newline
$b)$ One has $w_p(0) = w_p(1) = 1$ and $w_p(1-\alpha) = w_p(\alpha)~\forall \alpha\in I_p$.\newline
$c)$ For any fixed $n$, there are only finitely many $\alpha$'s for which $w_p(\alpha)\neq 0$ modulo $T^{n+1}\F_p[[T]]$.
\end{lem}
\proof $a)$ For each $n$ there exists at most one $k$ such that $p^n\alpha = k\in\N$. Such $k$ exists iff $n\ge n_0$ where $p^{n_0}$ is the denominator of the reduced form of $\alpha$. Thus the series \eqref{21m} is convergent in $\F_p[[T]]$.\newline
$b)$ One has $S_{p^n}(x,y) = S_{p^n}(y,x)$ so that $a(n,p^n-k) = a(n,k)$ and $w_p(1-\alpha)= w_p(\alpha)$. Moreover \eqref{1m} for $y_0=0$  gives the vanishing of $S_n(x,0)$ for all $n\ge 1$. Thus one gets that $a(n,0) = a(n,p^n)=0$ for all $n\ge 1$ and $w_p(0)=w_p(1)=1$.\newline
$c)$ If $w_p(\alpha)\neq 0$ mod $T^{n+1}\F_p[[T]]$, then $p^n\alpha\in\Z$ and thus $p^n\alpha\in\Z\cap[0,p^n]$ which is a finite set.
\endproof
We now consider the ring $S[[T]]$ of formal series in $T$ with coefficients in $S$. Since $S$ is of characteristic $p$, it contains $\F_p$ and one has
\begin{equation}\label{22m}
\F_p[[T]]\subset S[[T]]
\end{equation}
We use $\tau$ to obtain a bijection of $S[[T]]$ with $R$,
\begin{equation}\label{23m}
\tilde\tau(\sum s_nT^n) = \sum\tau(s_n)p^n\in R
\end{equation}
where the sum on the right is $p$-adically convergent.
\begin{thm}\label{teichthm} For any $x,y\in S$ one has
\begin{equation}\label{24m}
\tau(x)+\tau(y) = \tilde\tau(\sum_{\alpha\in I_p}w_p(\alpha)x^\alpha y^{1-\alpha})
\end{equation}
\end{thm}
\proof First, modulo $T^n$ the sum in the rhs is finite, thus
\begin{equation}\label{25m}
s=\sum_{\alpha\in I_p}w_p(\alpha)x^\alpha y^{1-\alpha}
\end{equation}
is convergent and defines an element of $S[[T]]$. By \eqref{21m} one has
\[
s = \sum_{\alpha\in I_p}\sum_n a(n,k)T^nx^\alpha y^{1-\alpha} = \sum_n T^n\sum_k a(n,k) x^{k/p^n} y^{1-k/p^n}
\]
which by \eqref{19m} gives $s = \sum S_{p^n}(x^{1/p^n}, y^{1/p^n})T^n$. Thus
\begin{equation}\label{26m}
\tilde\tau(s) = \sum\tau(S_{p^n}(x^{1/p^n},y^{1/p^n}))p^n
\end{equation}
and \eqref{13m} gives the required equality \eqref{24m}.
\endproof
 Note that \eqref{24m} extends to arbitrary sums $\sum \tau(x_i)$ for $x_i\in S$. One extends the definition of the $S_n(x,y)$ to $k$ variables by
 \begin{equation}\label{multipro}
\prod_1^k(1- tx_j)=\prod_{n=1}^\infty(1-S_n(x_1,\ldots,x_k)t^n)
\end{equation}
and one expands this polynomial for prime powers as
\begin{equation}\label{19multi}
S_{p^n}(x_1,\ldots,x_k) = \sum_{\sum m_j=p^n}a(n,m_1,\ldots, m_k)\prod x_j^{m_j}
\end{equation}
We let
$$
\bar a(n,m_1,\ldots, m_k) \equiv a(n,m_1,\ldots, m_k)~\text{mod}~p
$$
Then one considers the simplex
\begin{equation}\label{simplexp}
   I_p^{(k)}=\{(\alpha_j)\mid \alpha_j\in I_p,~\sum_1^k\alpha_j = 1\}
\end{equation}
and one defines a map from $I_p^{(k)}$ to $\F_p[[T]]$ by
\begin{equation}\label{21multi}
w_p: I_p^{(k)}\to \F_p[[T]],\quad w_p(\alpha_1,\ldots,\alpha_k)=\sum_{m_j/p^n=\alpha_j}\bar a(n,m_1,\ldots, m_k)T^n
\end{equation}

 \begin{thm}\label{teichthm1} For any $x_j\in S$ one has
\begin{equation}\label{32m}
\sum_{i=1}^k\tau(x_i) = \tilde\tau(\sum_{I_p^{(k)}}w_p(\alpha_1,\ldots,\alpha_k)\,x_1^{\alpha_1}\cdots x_k^{\alpha_k})
\end{equation}
\end{thm}

The proof is the same as for Theorem \ref{teichthm}.
Note that the coefficients $w_p(\alpha_1,\ldots,\alpha_k)$ {\em depend on $T$} and thus one is forced to extend the perfect ring $S$ in order to deform its operations. Evaluation at $T=0$, gives the usual rules of addition in $S$.

\section{Characteristic one}\label{charaone}

We refer to \cite{Golan} for the  general theory of semi-rings. In a semi-ring $R$ the additive structure $(R,+)$ is no longer that of an abelian group but is a commutative monoid with neutral element $0$. The multiplicative structure $(R,\cdot)$ is a monoid with identity $1\neq 0$ and distributivity holds while $r\cdot 0=0\cdot r =0$ for all $r\in R$.
A semi-ring  $R$ is called a {\em semi-field} when every non-zero element in $R$ has a multiplicative inverse, or equivalently when the set of non-zero elements in $R$ is a (commutative) group for the multiplicative law.\vspace{.05in}

\subsection{Perfect semi-rings of characteristic one}
 We let $\B=\{0,1\}$ be the only finite semi-field which is not a field. One has $1+1=1$ in $\B$ (\cite{Golan}, \cite{Lescot}).
 \begin{defn} \label{charisone} A semi-ring  $R$ is said to have characteristic $1$ when
\begin{equation}\label{idemplus}
1+ 1=1
\end{equation}
in $R$ \ie when $R$ contains $\B$ as the prime sub-semi-ring.
\end{defn}
It follows from distributivity  that one then has
\begin{equation}\label{27n}
a+a=a\qquad\forall a\in R
\end{equation}
 This justifies the term ``additively idempotent" frequently used in semi-ring theory as a synonymous for ``characteristic one".
  Note first that \eqref{27n} implies that
\begin{equation}\label{28n}
x+y = 0~\Rightarrow~ x=0,~y=0
\end{equation}
Indeed one has $x = x+(x+y) = (x+x)+y = x+y = 0$. \vspace{.05in}

A semi-ring $R$ of characteristic one inherits a canonical  partial order (\cf \cite{Golan}).
\begin{lem}\label{seclem} For $x$ and $y\in R$ one has
\begin{equation}\label{37n}
y+x = y~\Leftrightarrow~\exists z\in R,~y = x+z
\end{equation}
Moreover this relation defines a partial order $x\le y$ on $R$ and
\begin{itemize}
    \item $x_j\leq y_j$, $\forall j$, implies $\sum x_j\leq \sum y_j$.
  \item $x\leq y$ implies $ax\leq ay$ for all $a\in R$.
  \item For any $x,y\in R$ the sum $x+y$ is the least upper bound of $x$ and $y$.
\end{itemize}
\end{lem}
\proof If $y=x+z$, then $x+y = (x+x)+z = x+z = y$. Thus the equivalence \eqref{37n} follows.  The relation $x\le y$ iff $x+y=y$ is reflexive, transitive (if $x\leq y$ and $y\leq z$ then $x+z=x+(y+z) = (x+y)+z = y+z = z$ so that $x\leq z$) and anti-symmetric ($x\le y,~y\le x~\Rightarrow~x = x+y=y$). The compatibility with addition and multiplication is straightforward and finally if $x\leq z$ and $y\leq z$ one gets $x+y\leq z+z=z$ so the last property follows.
\endproof

By definition, a semi-ring $R$ is {\em multiplicatively cancellative}  when
\begin{equation}\label{29n}
x\neq 0,~xy = xz~\Rightarrow~ y=z
\end{equation}
Note that this condition means that the natural morphism
\begin{equation}\label{30n}
R \to S^{-1}R,\quad S = R\setminus\{0\}
\end{equation}
is an injection. Thus $R$ embeds in a semi-field of characteristic one. This set-up is thus closely related to that of a partially ordered group $G$.
We recall, from \cite{Golan} Propositions 4.43 and 4.44   the following result which describes the analogue of the Frobenius endomorphism in characteristic $p$.

\begin{lem}\label{theprop} Let $R$ be a \mc  semi-ring of characteristic one. Then the map $x \to x^n$ is an injective endomorphism of $R$ for any $n\in\N$.
\end{lem}
\proof Note first that for any $m$ one has
\begin{equation}\label{32n}
(x+y)^m = \sum_{k=0}^m x^k y^{m-k}
\end{equation}
as follows from \eqref{27n}.
To show that $x^n+y^n = (x+y)^n$ we first show
\begin{equation}\label{31n}
(x^n+y^n)(x+y)^{n-1} = (x+y)^{2n-1}
\end{equation}
which follows from \eqref{32n}  by noticing that if $x^a y^b$ is a monomial with $a+b = 2n-1$ then $a\ge n$ or $b\ge n$. From \eqref{31n} and the cancellative property one gets
\begin{equation}\label{33n}
(x+y)^n = x^n+y^n\qquad\forall x,y\in R
\end{equation}
Let us show that $x^n=y^n~\Rightarrow~x=y$. We start with $n=2$.

$x^2=y^2~\Rightarrow~x^2=x^2+y^2 = (x+y)^2 = x^2+xy+y^2 = x^2+xy$. Thus $xx = x(x+y)$ and canceling $x$ one gets $x = x+y$. Similarly $y = (x+y)$ and thus $y=x$. In general one has, assuming $x^n=y^n$
\begin{align*}
x^n &= (x+y)^n = x^n + yx^{n-1}+\cdots+y^{n-1}x+y^n = x^n + x^{n-1}y+\cdots xy^{n-1} =\\ &=x(x^{n-1}+x^{n-2}y+\cdots+y^{n-1}) = x(x+y)^{n-1}
\end{align*}
Thus since $R$ is \mc one gets $x^{n-1}=(x+y)^{n-1}$ and by induction  this gives $x=x+y$ and similarly $y = x+y$.\endproof

\begin{defn}\label{perfect}  Let $R$ be a \mc  semi-ring of characteristic one.
Then $R$ is {\em perfect}
 when for any $n$ the map $x\to x^n$ is surjective.
 \end{defn}

\begin{prop} \label{propfrob} Let $R$ be a \mc perfect semi-ring of characteristic one.
 The map $x\to x^n$  defines an automorphism $\theta_n\in\Aut(R)$ and
the equality
\begin{equation}\label{34n}
\theta_\alpha = \theta_a\theta_b^{-1},\quad \alpha=a/b\in\Q^*_+
\end{equation}
defines  an action of $\Q^*_+$ on $R$ fulfilling
\begin{equation}\label{all}
\theta_{\lambda\lambda'} = \theta_\lambda\circ\theta_{\lambda'},\quad \theta_\lambda(x)\theta_{\lambda'}(x) = \theta_{\lambda+\lambda'}(x)
\end{equation}
\end{prop}
\proof The first statement follows from Lemma \ref{theprop}. For $\alpha = a/b$, $\theta_\alpha(x)$ is the unique solution of
\[
z^b = x^a,\quad z\in R
\]
One has $\theta_{a_1}\theta_{a_2} = \theta_{a_1a_2}$ and the first part of \eqref{all} follows. For the second, note that if
\[
z_j^{b_j} = x^{a_j}, \quad j=1,2
\]
one gets $(z_1z_2)^{b_1b_2} = x^{a_1b_2+a_2b_1}$.
\endproof

In the above context we shall use the notation
\begin{equation}\label{xalpha}
  \theta_\alpha(x)=x^\alpha\qqq \alpha\in\Q^*_+, \ x\in R
\end{equation}

\begin{lem} \label{ineqpow} Let $R$ be a \mc perfect semi-ring of characteristic one. For $\alpha\in \Q\cap (0,1),~x,y\in R$ one has $x^\alpha y^{1-\alpha} \le x+y$.
\end{lem}
\proof  Let  $\alpha = \frac a n,~1-\alpha = \frac b n$ with $a+b = n$. Since $\theta_n$ is an automorphism,
\begin{equation}\label{38n}
x\le y~\Leftrightarrow~ x^n \le y^n\,.
\end{equation}
 Thus we just need to show that $x^a y^b \le (x+y)^n$ which follows from \eqref{32n}
\endproof

\begin{example}\label{mainex}{\rm We let $X$ be a compact space and $R= \tilde C(X,(0,\infty))$ be the semi-ring obtained by adjoining $0$ to the set $C(X,(0,\infty))$ of continuous functions  from $X$ to $(0,\infty)$ endowed with the operations
\begin{equation}\label{41n}
 (f+g)(x) = \text{Sup}(f(x),g(x)),~(fg)(x) = f(x)g(x)\qqq x\in X
\end{equation}
For each $x\in X$ the evaluation at $x$ gives a homomorphism from $R$ to $\rmax$, where $\rmax=\R_+=[0,\infty)$ endowed with the ordinary product and a new idempotent addition $x + y=\max\{ x ,y \}$. Note that it is not true in this example that an increasing bounded sequence $f_n$ will have a least upper bound since in general the latter will be given by a semi-continuous function. What is true however is that the intervals of the form
\begin{equation}\label{interval}
    [\rho_1,\rho_2]=\{x\in R\mid \rho_1\leq x\leq \rho_2\}\,, \ \rho_1\neq 0
\end{equation}
are {\em complete} for a suitable metric. We shall see in \S \ref{wittonesect} how to obtain a natural analogue of the $p$-adic metric in the context of characteristic one.
}
\end{example}

\subsection{Symmetrization}\label{symmsect}

In this section we recall a well known operation, called symmetrization, which associates a ring $A^\Delta$ to a semi-ring (\cf \cite{Golan}). Note that this operation always gives a trivial result when $A$ is of characteristic one and will only be used below in the deformed semi-rings.
Given a semi-ring $A$, its symmetrization $A^\Delta$ is obtained as follows:
\begin{prop} Let $A$ be a semi-ring.
\begin{enumerate}
\item The product $A\times A$ with the following operations is a semi-ring
\begin{equation}\label{65n}
(a,b)+(c,d) = (a+c,b+d),\quad (a,b)(c,d) = (ac+bd,ad+bc)
\end{equation}
\item The following subset $J$ is an ideal of $A\times A$
\begin{equation}\label{66n}
J = \{(a,a)|a\in A\}
\end{equation}
\item The quotient $(A\times A)/J$ is a ring, denoted $A^\Delta$.
\end{enumerate}
\end{prop}
\proof $(1)$~ The rules \eqref{65n} are those  of the group semi-ring $A[\Z/2\Z]$ generated by $A$ and $U$ commuting with $A$ and fulfilling $U^2=1$.

$(2)$~One has $J+J\subset J$. Let us check that $JA\subset J$. One has
\[
(a,a)(c,d)=(ac+ad,ad+ac)\in J
\]
$(3)$~The quotient of a semi-ring by an ideal $J$ is defined by the equivalence relation
\begin{equation}\label{67n}
\alpha~\sim~\beta~\Leftrightarrow~\exists j,j'\in J,~\alpha+j = \beta+ j'
\end{equation}
It is true in general that the quotient by an ideal is still a semi-ring \cite{Golan}. Let us show that the quotient $(A\times A)/J$ is a ring. To see this it is enough to show that any element $(a,b)$ has an additive inverse. But
\[
(a,b)+(b,a) = (a+b,a+b)\in J
\]
so that $(b,a)$ is the additive inverse of $(a,b)$.
\endproof
In the above case the equivalence relation \eqref{67n} means
\begin{equation}\label{68n}
(a,b)~\sim~(a',b')~\Leftrightarrow~\exists u,v\in A, \  a+u = a'+v,~b+u = b'+v
\end{equation}
This equivalence relation is the same as the one defining the Grothendieck group of the additive monoid which is given  by
\begin{equation}\label{69n}
(a,b)~\sim~(a',b')~\Leftrightarrow~\exists c,~a+b'+c = a'+b+c
\end{equation}
One checks that \eqref{69n} implies \eqref{68n}, taking $u = b'+c, v = b+c$. Conversely \eqref{68n} implies \eqref{69n}, taking $c = u+v$.

\section{Entropy and the $w(\alpha)$}\label{wittone}

Let $R$ be a \mc perfect semi-ring of characteristic one, and $\theta_\lambda(x) = x^\lambda$
be the automorphisms $\theta_\lambda\in \Aut(R)$ given by Proposition \ref{propfrob}.
We  let $I=\Q\cap(0,1)$ and consider a map
\begin{equation}\label{mapw}
    w:I\to R^\times\,, \ w(\alpha)\in R^\times \qqq \alpha\in I
\end{equation}
We extend $w$ to $\bar I=\Q\cap [0,1]$ by setting $w(0)=w(1)=1$.
We consider a sum of the form
\begin{equation}\label{1n}
x+_wy=\sum_{\alpha\in \bar I} w(\alpha)x^\alpha y^{1-\alpha}
\end{equation}
where, by convention we let $x^\alpha y^{1-\alpha}=x$ for $\alpha=1$ and $x^\alpha y^{1-\alpha}=y$ for $\alpha=0$.
 We first look formally at the associativity of the operation \eqref{1n}, disregarding the fact that it is an infinite sum. The commutativity means that
\begin{equation}\label{2n}
w(1-\alpha) = w(\alpha)\qqq \alpha\in I.
\end{equation}
One has:
\begin{align}\label{assocproof}
(x+_wy)+_wz &= (x+_wy) + z+ \sum_{\alpha\in I} w(\alpha)(x+_wy)^\alpha z^{1-\alpha} \\
 = x+y+z & +\sum_{\beta\in I} w(\beta)x^\beta y^{1-\beta}+\sum_{\alpha\in I} w(\alpha)\left(x+y+\sum_{\beta\in I} w(\beta)x^\beta y^{1-\beta} \right)^\alpha z^{1-\alpha} \nonumber \\
=x+y+z+T+&\sum_{\alpha, \beta\in I} w(\alpha)w(\beta)^\alpha x^{\alpha\beta} y^{\alpha(1-\beta)} z^{1-\alpha}\nonumber
\end{align}
where $T$ is the symmetric term given by
$$
T=\sum_{\beta\in I} w(\beta)x^\beta y^{1-\beta}+\sum_{\alpha\in I} w(\alpha)x^\alpha z^{1-\alpha}+\sum_{\alpha\in I} w(\alpha)y^\alpha z^{1-\alpha}
$$
By a similar computation one has
\[
x+_w(y+_wz) = x+y+z+T+\sum_{u,v\in I} w(u)w(v)^{1-u} x^u y^{v(1-u)} z^{(1-v)(1-u)}
\]
Thus it is enough to compare the two expressions
$$
\sum_{\alpha, \beta\in I} w(\alpha)w(\beta)^\alpha x^{\alpha\beta} y^{\alpha(1-\beta)} z^{1-\alpha}\,; \ \ \sum_{u,v\in I} w(u)w(v)^{1-u} x^u y^{v(1-u)} z^{(1-v)(1-u)}.
$$
Equating the exponents of $x,y,z$ one gets the transformation $u = \alpha\beta$, $v = \frac{\alpha(1-\beta)}{1-\alpha\beta}$. This transformation gives a bijection of $I\times I$ to $I\times I$ whose inverse is given by
$$
\alpha=u+v-uv\,, \ \ \beta =\frac{u}{u+v-uv}.
$$
Thus the associativity holds provided
\begin{equation}\label{3n}
w(\alpha)w(\beta)^\alpha = w(\alpha\beta)w\left(\frac{\alpha(1-\beta)}{1-\alpha\beta}\right)^{1-\alpha\beta}\qqq \alpha,\beta\in I.
\end{equation}

\subsection{General solutions}
Both \eqref{2n} and \eqref{3n}  only involve the multiplicative group $G=R^\times$ of the semi-ring $R$ and we thus fix a multiplicative group $G$ and assume that it is uniquely divisible \ie that  $x^\alpha\in G$   makes sense for any $x\in G$ and positive rational number $\alpha$.
We now look for solutions of \eqref{2n} and \eqref{3n} in this framework. For each $n$ we use the notation \eqref{simplex}
\begin{equation}\label{simplexbis}
    \Sigma_n=\{(\alpha_1,\ldots,\alpha_n)\in I^{n}\mid \sum\alpha_j = 1\}
\end{equation}

\begin{lem} \label{prepare} Let $G$ be a uniquely divisible (multiplicative) group and $w(\alpha)\in G$, $\alpha\in I$ satisfy \eqref{2n} and \eqref{3n}, then\vspace{.05in}

$(1)$~ The following map from $\Sigma_n$ to $G$ is symmetric
\begin{align}\label{4n}
w(\alpha) &= w(\alpha_1)w\left(\alpha_2/(1-\alpha_1)\right)^{1-\alpha_1}
w\left(\alpha_3/(1-\alpha_1-\alpha_2)\right)^{1-\alpha_1-\alpha_2}\cdots\\
&\cdots w\left(\alpha_{n-1}/(1-\alpha_1-\cdots -\alpha_{n-2})\right)^{1-\alpha_1-\cdots-\alpha_{n-2}}\qqq \alpha=(\alpha_1,\ldots,\alpha_n)\in \Sigma_n.\nonumber
\end{align}
$(2)$~ Let $(J_k)$, $k=1,\ldots m$, be a partition of $\{1,\ldots ,n\}=\cup J_k$ then  one has
\begin{equation}\label{5n}
w(\alpha) = w(\beta)\prod_1^m w(\gamma_k)^{\beta_k}\qqq \alpha \in \Sigma_n
\end{equation}
where $\beta_k=\sum_{J_k}\alpha_j$ and\footnote{when $J_k$ has only one element $w(\gamma_k)=1$} $(\gamma_k)_j=\alpha_j/\beta_k$ for all $j\in J_k$.

$(3)$~There exists a unique homomorphism $\chi:\Q_+^*\to G$, such that
\begin{equation}\label{6n}
\chi(n^{-1}) = w(n^{-1},\cdots,n^{-1})\qqq n\in \N.
\end{equation}
$(4)$~ One has
\begin{equation}\label{7n}
w(\beta) = \chi(\beta)^\beta\chi(1-\beta)^{1-\beta}\qqq \beta \in I
\end{equation}
\end{lem}
\proof $(1)$~ Permuting $\alpha_1$ and $\alpha_2$ only affects the first two terms. Let us show  that
\begin{equation}\label{transpo}
w(\alpha_1)w(\alpha_2/(1-\alpha_1))^{1-\alpha_1} = w(\alpha_2)w(\alpha_1/(1-\alpha_2))^{1-\alpha_2}
\end{equation}
 Using \eqref{3n} for $\alpha = 1-\alpha_1,~\beta = \alpha_2/(1-\alpha_1)$ one gets $$\alpha\beta = \alpha_2\,, \  1-\frac{\alpha(1-\beta)}{1-\alpha\beta} = \frac{1-\alpha}{1-\alpha\beta} = \frac{\alpha_1}{1-\alpha_2}$$
 Thus \eqref{transpo} follows from \eqref{2n} and \eqref{3n}.
In the same way let us show that, for $k<n-1$, the permutation $\alpha_k~\leftrightarrow~\alpha_{k+1}$ does not change \eqref{4n}. It only  affects the following two consecutive terms
\[
T=w(\alpha_k/(1-(\alpha_1+\ldots+\alpha_{k-1})))^{1-\alpha_1-\cdots-\alpha_{k-1}}
w(\alpha_{k+1}/(1-(\alpha_1+\ldots+\alpha_{k})))^{1-\alpha_1-\cdots-\alpha_{k}}
\]
Moreover one has
$$
T=(w(a)w(b/(1-a))^{1-a})^\beta
$$
where
$$
\beta=1-\alpha_1-\cdots-\alpha_{k-1}\,, \ a=\alpha_k/(1-(\alpha_1+\ldots+\alpha_{k-1}))\,,
b =\alpha_{k+1}/(1-(\alpha_1+\ldots+\alpha_{k-1}))
$$
since
$$
b/(1-a)=\alpha_{k+1}/(1-(\alpha_1+\ldots+\alpha_{k}))\,.
$$
Thus the permutation $\alpha_k~\leftrightarrow~\alpha_{k+1}$ interchanges $a~\leftrightarrow~b$ and the invariance follows from \eqref{transpo}. Finally the invariance under the permutation $\alpha_{n-1}~\leftrightarrow~\alpha_n$ follows from \eqref{2n}.

$(2)$~  We first consider the special case of the partition $\{1,\ldots,n\}=\{1,2\}\cup_k \{k\}$. Then  \eqref{5n} takes the form
\begin{equation}\label{special}
    w(\alpha_{1},\alpha_{2},\ldots,\alpha_n)=w(\alpha_1+\alpha_2,\alpha_3,
    \ldots,\alpha_n)w(\alpha_{1}/(\alpha_1+\alpha_2))^{\alpha_1+\alpha_2}
\end{equation}
In that case using \eqref{4n}, equation \eqref{5n} reduces to
 $$
w(u)w(v/(1-u))^{1-u} =w(u+v) w(u/(u+v))^{u+v}
$$
for $u=\alpha_{1}$, $v=\alpha_{2}$ (so that $u+v\leq 1$).
In turns this follows from \eqref{3n}, for $\alpha=u+v$, $\beta=u/(u+v)$.

 To prove \eqref{5n} in general, one proceeds by induction on $n$. If all the $J_k$ have only one element the equality is trivial. Let then $J_k$ have at least two elements, which using the symmetry proved in $(1)$ can be assumed to be $1,2\in \{1,\ldots,n\}$. One shows that the replacement $(\alpha_{1},\alpha_{2})\mapsto \alpha_{1}+\alpha_{2}$ in both sides of \eqref{5n} has the same effect. Using \eqref{special} one checks that both sides are divided by the factor $$w(\alpha_{1}/(\alpha_{1}+\alpha_{2}))^{\alpha_{1}+\alpha_{2}}$$
 Indeed for the left hand side this follows from \eqref{special} and in the right hand side only the term $w(\gamma_k)^{\beta_k}$ gets modified to $w(\gamma'_k)^{\beta_k}$ and one has, using \eqref{special},
 $$
 w(\gamma'_k)=w(\gamma_k)w(\alpha_{1}/(\alpha_{1}+\alpha_{2}))^{-(\alpha_{1}+\alpha_{2})/\beta_k}\,.
 $$
This proves \eqref{5n} by induction on $n$ since it holds for $n=2$ as one checks directly for the two possible partitions.

$(3)$~Let $\alpha \in \Sigma_n$ and $\alpha'\in \Sigma_m$. Then the $\alpha_{ij} = \alpha_i\alpha_j'$ belong to $I$ and add up to $1$. Then all the $\gamma_k$ are the same and \eqref{5n} gives
\begin{equation}\label{8n}
w((\alpha_i\alpha_j')) = w((\alpha_i))w((\alpha_j'))
\end{equation}
One applies this for $\alpha_i = 1/n$, $i=1,\ldots, n$ and $\alpha_j' = 1/m$, $j=1,\ldots, m$ and gets using \eqref{6n},
\begin{equation}\label{9n}
\chi((nm)^{-1}) = \chi(n^{-1})\chi(m^{-1})\qqq n,m \in \N
\end{equation}
Thus there exists a unique extension of $\chi$ to a homomorphism from $\Q_+^*$, which is determined by its value on $1/p$, $p$ a prime number.

$(4)$~ Let $\beta = \frac nm$. Take $\alpha\in \Sigma_{m}$ with $\alpha_j = \frac 1m$ for all $j$, and the partition in two sets $J_1$ and $J_2$ containing respectively $n$ and $m-n$ elements. Then apply \eqref{5n}. The left hand side is $\chi(\frac 1m)$ by construction. The first term in the right hand side is $w(\frac nm)=w(\beta)$. The term $w(\gamma_1)$ gives $\chi(\frac 1n)$ and similarly $w(\gamma_2)$ gives $\chi(\frac{1}{m-n})$.  Thus \eqref{5n} gives $\chi(\frac 1m) = w(\beta)\chi(\frac 1n)^\beta\chi(\frac{1}{m-n})^{1-\beta}$ and thus $w(\beta) = \chi(\frac{n}{m})^\beta\chi(\frac{m-n}{m})^{1-\beta}= \chi(\beta)^\beta\chi(1-\beta)^{1-\beta}$.
\endproof

\begin{thm} \label{checkchar} Let $\chi$ be a homomorphism $\Q_+^\times \to G$, then the function
\begin{equation}\label{10n}
w(\alpha) = \chi(\alpha)^\alpha\chi(1-\alpha)^{1-\alpha}
\end{equation}
fulfills \eqref{2n} and \eqref{3n}. Moreover all solutions of \eqref{2n} and \eqref{3n} are of this form.
\end{thm}
\proof Let us first show that \eqref{10n} implies \eqref{3n}. One has
\[
w(\alpha\beta)w(\frac{\alpha(1-\beta)}{1-\alpha\beta})^{1-\alpha\beta} = \chi(\alpha\beta)^{\alpha\beta}\chi(1-\alpha\beta)^{1-\alpha\beta}
\]
\[
\times\ \ \chi\left(\frac{\alpha(1-\beta)}{1-\alpha\beta}\right)^{\alpha(1-\beta)}
\chi\left(\frac{1-\alpha}{1-\alpha\beta}\right)^{1-\alpha}
\]
which using  the multiplicativity of $\chi$ gives
\[
\chi(\alpha)^{\alpha\beta}\chi(\beta)^{\alpha\beta}
\chi(1-\alpha\beta)^{(1-\alpha\beta-\alpha(1-\beta)-(1-\alpha))}\chi(\alpha)^{\alpha(1-\beta)}\chi(1-\beta)^{\alpha(1-\beta)}
\]
\[
\times\ \ \chi(1-\alpha)^{1-\alpha}=\chi(\alpha)^\alpha\chi(1-\alpha)^{1-\alpha}\chi(\beta)^{\alpha\beta}\chi(1-\beta)^{\alpha(1-\beta)}
\]
which agrees with
\[
w(\alpha)w(\beta)^\alpha = \chi(\alpha)^\alpha\chi(1-\alpha)^{1-\alpha}\chi(\beta)^{\alpha\beta}
\chi(1-\beta)^{\alpha(1-\beta)}\,.
\]
Conversely by Lemma \ref{prepare} $(4)$~ all solutions of \eqref{2n} and \eqref{3n} are of the form
\eqref{10n}.
\endproof

\subsection{Positive solutions}
We let, as above, $R$ be a \mc perfect semi-ring of characteristic one.
The uniquely divisible group $G=R^\times$ is a vector space over $\Q$ using the action $\theta_\alpha(x)=x^\alpha$ and is partially ordered by Lemma \ref{seclem}. We shall now make the stronger assumption that it is a partially ordered vector space over $\R$. Thus $G$ is a partially ordered group  endowed with a one parameter group of automorphisms $\theta_\lambda\in\text{Aut}(G)$, $\lambda\in\R^\times$ such that, with $\theta_0(x)=1$ for all $x$ by convention,
\begin{equation}\label{24n}
\theta_{\lambda\lambda'} = \theta_\lambda\circ\theta_{\lambda'},\quad \theta_\lambda(x)\theta_{\lambda'}(x) = \theta_{\lambda+\lambda'}(x)
\end{equation}
We assume the following compatibility (closedness) of the partial order with the vector space structure
\begin{equation}\label{26n}
\lambda_n \to \lambda,\quad \theta_{\lambda_n}(x)\ge y ~\Rightarrow~\theta_\lambda(x)\ge y
\end{equation}

\begin{thm}\label{entropythm} Let $w: I\to G$ fulfill \eqref{2n} and \eqref{3n} and
\begin{equation}\label{11n}
w(\alpha)\geq 1,\qquad\forall\alpha\in I\,.
\end{equation}
Then there exists $\rho\in G$, $\rho\ge 1$ such that
\begin{equation}\label{23n}
w(\alpha)=\rho^{ S(\alpha)},~S(\alpha) = -\alpha\log(\alpha)-(1-\alpha)\log(1-\alpha)\qqq\alpha\in I.
\end{equation}
\end{thm}

\proof
We use additive notations so that for $x\in G$, one has $\log(x)\in E$ where $E$ is a partially ordered vector space over $\R$. We let  $s(\alpha)=\log(w(\alpha))\ge 0$ and $L(\alpha) = \log(\chi(\alpha))$ for $\alpha\in \Q_+^*$. Note that by \eqref{4n} and \eqref{11n} one gets
\begin{equation}\label{12n}
w(\alpha_1,\ldots,\alpha_n)\ge 1
\end{equation}
and thus by \eqref{6n}, $\chi(1/n)\ge 1$, so that
\begin{equation}\label{13n}
L(1/n)\ge 0\qqq n\in \N
\end{equation}
One lets $l(p) = L(1/p)\ge 0$, for each prime $p$.
By \eqref{7n} one has
\begin{equation}\label{14n}
s(\alpha) = \alpha L(\alpha)+(1-\alpha)L(1-\alpha),\qquad  s(\alpha)\geq 0\qqq \alpha\in I.
\end{equation}
Let $p_1$ and $p_2$ be two primes and $n_j$ integers such that $\alpha = \frac{p_1^{n_1}}{p_2^{n_2}}<1$. One has: $1-\alpha = \frac{p_2^{n_2}-p_1^{n_1}}{p_2^{n_2}}$ and if one lets
\begin{equation}\label{15n}
p_2^{n_2}-p_1^{n_1} = \prod q_j^{a_j}
\end{equation}
be the prime factor decomposition of $p_2^{n_2}-p_1^{n_1}>0$, one gets
\begin{equation}\label{16n}
L(1-\alpha) = n_2l(p_2)-\sum a_j l(q_j)\le n_2 l(p_2)
\end{equation}
Thus, since $s(\alpha)\geq 0$, $\alpha L(\alpha) = s(\alpha)-(1-\alpha)L(1-\alpha)$ fulfills
\begin{equation}\label{17n}
\alpha L(\alpha)\ge -(1-\alpha)L(1-\alpha)\ge -(1-\alpha)n_2 l(p_2)
\end{equation}
But $L(\alpha) = n_2 l(p_2)-n_1 l(p_1)$ and thus dividing by $n_2$ we get
\begin{equation}\label{18n}
\alpha(l(p_2)-\frac{n_1}{n_2}l(p_1))\ge -(1-\alpha)l(p_2)
\end{equation}
which gives
\begin{equation}\label{20n}
l(p_2) - \alpha\frac{n_1}{n_2} l(p_1)\ge 0\,.
\end{equation}
Let $$a=\frac{\log(p_2)}{\log(p_1)}$$ then $a\notin \Q$.
Taking rational approximations (using the density of $\Z+a\Z$ in $\R$), one gets a sequence $(n_1(j),n_2(j))$ such that
$$
n_1(j)-a\, n_2(j)<0\,, \ \ n_1(j)- a \, n_2(j)\to 0\,, \ \text{when} \ j\to \infty
$$
one then gets, when $j\to\infty$ and with $n_k=n_k(j)$,
\begin{equation}\label{ratappr}
 \frac{n_1}{n_2}\to a=\frac{\log(p_2)}{\log(p_1)}\,, \ \  \frac{p_1^{n_1}}{p_2^{n_2}}<1\,, \ \frac{p_1^{n_1}}{p_2^{n_2}}\to 1
\end{equation}
and using \eqref{26n} one gets from \eqref{20n},
\begin{equation}\label{21n}
\frac{l(p_2)}{\log(p_2)}-\frac{l(p_1)}{\log(p_1)}\ge 0\,.
\end{equation}
Exchanging the roles of $p_1$ and $p_2$ thus gives the equality
\begin{equation}\label{22n}
\frac{l(p_2)}{\log(p_2)}=\frac{l(p_1)}{\log(p_1)}\qqq p_1,p_2\,.
\end{equation}
This shows  that $\lambda=l(p)/\log(p)$ is a positive element of $E$ independent of the prime $p$. One then has
\[
L(1/n) = \lambda\log(n)\qqq n\in\N,
\]
and thus $L(\alpha) = -\lambda\log(\alpha)$ for $\alpha\in\Q^*_+$. One thus gets
\[
\log(w(\alpha)) = \lambda(-\alpha\log(\alpha)-(1-\alpha)\log(1-\alpha)) = \lambda S(\alpha)
\]
which gives \eqref{23n}, using the multiplicative notation.
\endproof

Note that the function $w$ automatically extends by continuity from $I$ to $[0,1]$.

\section{Analogue of the Witt construction in characteristic one}\label{wittonesect}

 Let $R$ be a \mc perfect semi-ring of characteristic one. We look at the meaning of an expression of the form
\begin{equation}\label{35n}
\sum_{\alpha\in \bar I}w(\alpha)x^\alpha y^{1-\alpha}
\end{equation}
Here, $I$ is $\Q\cap(0,1)$, $\bar I=\Q\cap [0,1]=I\cup \{0,1\}$ and we use the notation $x^\alpha = \theta_\alpha(x)$ (\cf \eqref{xalpha}). Note that this does not make sense for $\alpha \in \{0,1\}$. Thus we let
\begin{equation}\label{36n}
x^\alpha y^{1-\alpha} = x\quad\text{for $\alpha=1$}, \ x^\alpha y^{1-\alpha} =  y\quad\text{for $\alpha=0$}
\end{equation}
The index set $\bar I$ is countable and we treat \eqref{35n} as a discrete sum.
\begin{lem} \label{divis} Let $I(n) = \frac 1 n \Z\cap[0,1]$. Then the partial sums
\begin{equation}\label{39n}
s(n) = \sum_{I(n)}w(\alpha) x^\alpha y^{1-\alpha}
\end{equation}
form an increasing family for the partial order given by divisibility $n|m$.
\end{lem}
\proof We need to show that $s(n)\le s(m)$ when $n$ divides $m$. In that case one has $I(n)\subset I(m)$ and thus the conclusion follows from \eqref{37n}.
\endproof
Moreover, assuming $w(\alpha)\le \rho$ for all $\alpha\in I$ and some fixed $\rho\in R$, we get using  Lemmas~\ref{seclem} and \ref{ineqpow}  that
\begin{equation}\label{40n}
x+y\leq s(n) \le  (x+y)\rho \qquad\forall n.
\end{equation}
The hypothesis that any increasing bounded sequence has a least upper bound is too strong since it fails in example \ref{mainex} for instance. We shall now show  that the $s(n)$ form a Cauchy sequence and converge to the least upper bound of the $s(n)$ provided one passes to the completion for the
$\rho$-adic  metric which we now construct.

\subsection{Completion for the $\rho$-adic  metric}\label{completionsect}
 Let $R$ be a \mc perfect semi-ring of characteristic one. Let $\rho\in R^\times,~\rho\ge 1$.
We want to use the $w(\alpha)$ given by \eqref{23n} and for this we first need to perform a suitable completion. The function $S(\alpha)$ is positive and bounded by $\log 2<1$  for $\alpha\in I$. By \eqref{40n} we just need to complete the intervals  of the form
\begin{equation}\label{inter}
[\rho^{-n},\rho^n] = \{x\in R\mid \rho^{-n}\leq x\leq \rho^n \},\quad n\in \N
\end{equation}
We let
\begin{equation}\label{rrho}
    R_\rho=\{0\}\cup \bigcup_{n\in \N} [\rho^{-n},\rho^n]\subset R.
\end{equation}

\begin{lem}\label{rhoadd} $(1)$~For any $\alpha,\beta\in \Q_+^*$
\begin{equation}\label{56n}
\rho^\alpha + \rho^\beta = \rho^{\alpha\vee\beta},\quad \alpha\vee\beta = \text{Sup}(\alpha,\beta).
\end{equation}
$(2)$~$R_\rho$ is a perfect sub semi-ring of $R$.
\end{lem}
\proof
$(1)$~Let us show that $\rho^\alpha\leq \rho^\beta$ for $\alpha<\beta$. It is enough to show that $1+\rho^{\beta-\alpha} = \rho^{\beta-\alpha}$, \ie it is enough to show that
\begin{equation}\label{57n}
1+\rho^\gamma = \rho^\gamma,\quad\forall\gamma\in\Q^+
\end{equation}
For $\gamma = \frac k n$ this follows from the additivity of $\theta_\gamma$ since $1+\rho = \rho$.

$(2)$~It follows from $(1)$ that $R_\rho$ is a sub semi-ring of $R$. Moreover for $\alpha\in \Q_+^*$
and $x\in [\rho^{-n},\rho^n]$ one has $\theta_\alpha(x)\in [\rho^{-n\alpha},\rho^{n\alpha}]\subset R_\rho$ using $(1)$.
\endproof
We want to construct  a  metric $d$ on $R_\rho$ such that
\begin{equation}\label{hyp}
d(\rho^\alpha, 1)\to 0,\quad\text{when $\alpha\in \Q$},\quad \alpha\to 0.
\end{equation}
For $\alpha\in \Q_+^*$ we let
\begin{equation}\label{entourage}
    \cU_\alpha=\{(x,y)\mid x\leq y\rho^\alpha,\ y\leq x\rho^\alpha\}
\end{equation}
\begin{lem}\label{uniformstr} One has
\begin{equation}\label{unif}
    \cU_\alpha^{-1}=\cU_\alpha\,, \ \ \cU_\alpha\circ\cU_\beta\subset \cU_{\alpha+\beta}
\end{equation}
and
\begin{equation}\label{entsum}
    (x_j,y_j)\in \cU_\alpha\implies (\sum x_j\, ,\sum y_j)\in \cU_\alpha.
\end{equation}
\end{lem}
\proof For $(x,y)\in \cU_\alpha$ and $(y,z)\in \cU_\beta$ one gets
$$
x\leq y\rho^\alpha, \ y\leq z\rho^\beta \implies x\leq z\rho^{\alpha+\beta}
$$
and one gets in the same way $z\leq x\rho^{\alpha+\beta}$. This gives \eqref{unif}. Similarly
$$
x_j\leq y_j \rho^\alpha\implies \sum x_j\leq (\sum y_j)\rho^\alpha.
$$
which gives \eqref{entsum}.\endproof

We now define $d(x,y)\in \R_+$ for $x,y\in R_\rho\setminus \{0\}$ as follows
\begin{equation}\label{dxy}
    d(x,y)=\inf\{\alpha\mid (x,y)\in \cU_\alpha\}.
\end{equation}
Note that $d(x,y)<\infty$ for all $x,y\in R_\rho\setminus \{0\}$ since  if they both belong to some interval $[\rho^{-n},\rho^n]$ one has $d(x,y)\leq 2n$.

\begin{lem}\label{uniformdis} One has
\begin{equation}\label{triang}
    d(x,y)=d(y,x)\,, \ \ d(x,z)\leq d(x,y)+d(y,z)\qqq x,y,z\neq 0.
\end{equation}
\begin{equation}\label{addi}
   d(\sum x_i\, , \sum y_i)\le \text{Sup}(d(x_i,y_i))\qqq x_i,y_i\neq 0.
\end{equation}
\begin{equation}\label{multi}
    d(xy,zt)\leq d(x,z)+d(y,t)\qqq x,y,z,t \neq 0.
\end{equation}
\begin{equation}\label{power}
    d(x^\alpha,y^\alpha)\leq \alpha \,d(x,y)\qqq x,y\neq 0.
\end{equation}

\end{lem}
\proof \eqref{triang} follows from \eqref{unif}. Similarly \eqref{addi} follows from \eqref{entsum}.
Note that one has
$$
x\leq z\rho^\alpha\,, \ y\leq t\rho^\beta\implies xy\leq zt \rho^{\alpha+\beta}
$$
which implies \eqref{multi}. Finally
$$
x\leq y\rho^\beta\implies x^\alpha\leq y^\alpha\rho^{\alpha\beta}
$$
which implies \eqref{power}.
 \endproof
\begin{prop}\label{completion} The separated completion $\bar R_\rho$ obtained by adjoining $0$ to the separated completion of $(R_\rho\setminus \{0\},d)$ is a perfect semi-ring of characteristic one.

The action $\theta_\alpha$, $\alpha\in\Q_+^*$ on $\bar R_\rho$ extends by continuity to an action of $\R_+^*$ by automorphisms $\theta_\lambda\in \Aut(\bar R_\rho)$.
\end{prop}
\proof Lemma \ref{uniformdis} shows that $d$ defines a pseudo-metric on $R_\rho\setminus \{0\}$. In fact it also defines a pseudo-metric on the set of Cauchy sequences of $R_\rho\setminus \{0\}$. The quotient by the equivalence relation $d(x,y)=0$ is naturally endowed with addition using \eqref{addi} which implies
$$
d(x+y,z+t)\leq \text{Sup} (d(x,z),d(y,t))
$$
and multiplication using \eqref{multi} so that
$$
d(xy,zt)\leq d(x,z)+d(y,t)
$$
is still valid for Cauchy sequences. This shows that $\bar R_\rho$ is a semi-ring. The inequality \eqref{power} shows that the $\theta_\alpha$ extend to the completion (for $\alpha\in\Q_+^*$) by uniform continuity. The equality $1+1=1$ continues to hold in $\bar R_\rho$ which is thus a perfect semi-ring of characteristic one. Let us show that the action $\theta_\alpha$, $\alpha\in\Q_+^*$ on $\bar R_\rho$ extends by continuity to an action of $\R_+^*$.
First for any $x\in [\rho^{-n},\rho^n]$  and $\alpha\in\Q_+^*$ one has $d(x,1)\leq n$ and using \eqref{power}
\begin{equation}\label{unifor}
    d(x^\alpha,1)\leq n\,\alpha
\end{equation}
Using \eqref{multi} this implies
\begin{equation}\label{unifor1}
    d(x^\alpha,x^\beta)\leq n|\alpha-\beta|
\end{equation}
and shows that for $\alpha_j\to \lambda\in  \R_+^*$ the sequence $x^{\alpha_j}$ is a Cauchy sequence. \endproof
 We continue to use the notation $\theta_\lambda(x)=x^\lambda$ for $\lambda\in  \R_+^*$.

 \subsection{Construction of $W(R,\rho)$}\label{constructW}

We can now show that the $s(n)$ of Lemma \ref{divis} form a Cauchy sequence.
\begin{lem}\label{seclemcau} Let $R$, $\rho$ and $\bar R_\rho$ be as above.
Let $w(\alpha) = \rho^{S(\alpha)}$ for all $\alpha\in I$. Then for any $x,y\in \bar R_\rho$ the sequence
\begin{equation}\label{Cauchysequ}
s(n)=\sum_{I(n)}w(\alpha)x^\alpha y^{1-\alpha}
\end{equation}
is  a Cauchy sequence which converges to the lowest upper bound of the $w(\alpha)x^\alpha y^{1-\alpha}$.
\end{lem}
\proof We can assume $x,y\neq 0$ since otherwise $s(n)$ is constant. We estimate $d(s(n),s(nm))$. For this we write the elements of $I(nm)$ in the form
\begin{equation}\label{47n}
\alpha = \frac a n + \frac{k}{nm},\quad 0\le k\le m
\end{equation}
Using \eqref{addi} it is enough to estimate uniformly
\begin{equation}\label{48n}
d(\rho^{S(\alpha)}x^\alpha y^{1-\alpha}, \rho^{S(\frac a n)} x^{\frac a n} y^{1-\frac a n})
\end{equation}
Thus it is enough to show that for any $\epsilon >0$, there exists $\delta>0$ such that
\begin{equation}\label{52n}
d(\rho^{S(\alpha)}x^\alpha y^{1-\alpha}, \rho^{S(\beta)}x^\beta y^{1-\beta})\le \epsilon,~\forall\alpha,\beta\in[0,1],~|\alpha-\beta|<\delta
\end{equation}
This follows from \eqref{multi} which allows one to consider each of the three terms separately. One uses  \eqref{unifor1} for the terms $x^\alpha$, $ y^{1-\alpha}$ and  the uniform continuity of the function $S(\alpha)$, $\alpha\in[0,1]$ for the term $\rho^{S(\alpha)}$. By Lemma \ref{divis} the limit of the $s(n)$ gives the lowest upper bound of the $s(n)$ and hence of the $\rho^{S(\alpha)}x^\alpha y^{1-\alpha}$.
\endproof
Thus, under the hypothesis of Lemma \ref{seclemcau}, we get the convergence in $\bar R_\rho$. We keep a notation closely related to the Witt case and let
\begin{equation}\label{55n}
\sum_{\alpha\in \bar I}w(\alpha)x^\alpha y^{1-\alpha} = \lim_{n\to \infty} s(n)\,, \ \ s(n)= \sum_{\alpha\in I(n)} w(\alpha) x^\alpha y^{1-\alpha}
\end{equation}
Before we go any further we shall evaluate the new operation on the powers of $\rho$.
Let us evaluate  \eqref{55n} for $x = \rho^a,~y = \rho^b$. One has
\[
s(n)=\sum_{\alpha\in I(n)} \rho^{(S(\alpha)+\alpha a + (1-\alpha)b)}
\]
and thus, by \eqref{56n}, it is given by $\rho^{\sigma_n(a,b)}$ where
\begin{equation}\label{58n}
\sigma_n(a,b) = \text{Sup}_{\alpha\in I(n)}\left( S(\alpha) + \alpha a + (1-\alpha) b\right)
\end{equation}
\begin{lem}\label{tlem} For $a,b\in \R$, one has
\begin{equation}\label{59n}
\sigma_n(a,b)\to \log(e^a+e^b)\quad \text{when $n\to \infty$}
\end{equation}
\end{lem}
\proof The function $S(\alpha)+\alpha a + (1-\alpha)b = f(\alpha)$ is uniformly continuous for $\alpha \in[0,1]$ and thus $\sigma_n(a,b)$ tends to its maximum. One has
\begin{equation}\label{60n}
f'(\alpha) = \log(\frac{1-\alpha}{\alpha}) + a - b
\end{equation}
and $f''(\alpha) = - \frac{1}{\alpha(1-\alpha)}<0$, so that $f$ is strictly concave on $[0,1]$. Take $\alpha = \frac{e^a}{e^a+e^b}$. Then $1-\alpha = \frac{e^b}{e^a+e^b}$, $\frac{1-\alpha}{\alpha}= e^{(b-a)}$ and thus, by \eqref{60n}, $f'(\alpha) = 0$. Thus $f$ is strictly increasing in $[0,\alpha]$ and decreasing in $[\alpha,1]$. Its value at $\alpha$ is
\begin{align*}
&-\alpha\log(\alpha)-(1-\alpha)\log(1-\alpha)+\alpha a + (1-\alpha) b =\\&= -\alpha\log(e^a)-(1-\alpha)\log(e^b)+\alpha a+(1-\alpha)b + \alpha\log(e^a+e^b)+\\&+(1-\alpha)\log(e^a+e^b) = \log(e^a+e^b)
\end{align*}
\endproof

\begin{prop}\label{semiringexist} $(a)$~The following defines an associative operation in $\bar R_\rho$
\begin{equation}\label{63n}
x+_wy = \sum_{\alpha\in \bar I} w(\alpha)x^\alpha y^{1-\alpha},\quad 0+_wy = y+_w0 = y
\end{equation}
$(b)$~One has $(x+_wy)z = xz+_wyz$, $\forall x,y,z\in \bar R_\rho$.

$(c)$~One has
\begin{equation}\label{addrho}
    \rho^a+_w\rho^b=\rho^c\,, \ c=\log(e^a+e^b)\qqq a,b\in \R.
\end{equation}
\end{prop}
\proof By  Lemma~\ref{seclemcau} \eqref{63n} is well defined. The commutativity follows from the symmetry $w(1-\alpha)=w(\alpha)$. For $x=0$ all terms $x^\alpha y^{1-\alpha}$ vanish except when $\alpha=0$ which gives $y$, thus the sum gives $x+_wy =y$. To show associativity we proceed as in \eqref{assocproof} and we can now justify the equality
$$
\sum_{\alpha\in I}w(\alpha)\left(\sum_{\beta\in I} w(\beta)x^\beta y^{1-\beta} \right)^\alpha z^{1-\alpha}=
\sum_{\alpha, \beta\in I} w(\alpha)w(\beta)^\alpha x^{\alpha\beta} y^{\alpha(1-\beta)} z^{1-\alpha}
$$
using the continuity of the automorphisms $\theta_\alpha$ (\cf \eqref{power}).
 The equality $(b)$ follows from $z = z^\alpha z^{1-\alpha}$.
$(c)$~follows from Lemma  \ref{tlem}.
\endproof

\begin{cor}\label{semiringexist1}
Assume $\rho\neq 1$. Then $(\bar R_\rho, +_w, \cdot)$ is a semi-ring and the map
\begin{equation}\label{inclus}
    s\in \R_+\mapsto r(s)=\rho^{\log(s)}\,, \ \ r(0)=0,
\end{equation}
is an injective homomorphism of the semi-ring $\R_+$ (with ordinary addition) in
$(\bar R_\rho, +_w, \cdot)$.
\end{cor}
\proof Proposition \ref{semiringexist} shows that $(\bar R_\rho, +_w, \cdot)$ is a semi-ring. Using \eqref{addrho} it just remains to show that $\rho^a\neq \rho^b$ if $a\neq b$. One has $\rho^\alpha\neq 1$ for $\alpha\in \Q_+^*$, $\alpha\neq 0$ since $\rho\neq 1$ and $\theta_\alpha$ is an automorphism. By Lemma \ref{rhoadd}, this shows that $d(\rho^\alpha,\rho^\beta)=|\alpha-\beta|$ for $\alpha,\beta\in \Q$ and extends to
$$
d(\rho^\alpha,\rho^\beta)=|\alpha-\beta|\qqq \alpha,\beta\in \R
$$
which gives the injectivity of $r$.\endproof

We now define $W(R,\rho)$.
\begin{defn}\label{defnwr} Let  $R$ be a \mc perfect semi-ring of characteristic one and $\rho\in R$, $\rho\geq 1$. We let
\begin{equation}\label{70n}
W(R,\rho) = (\bar R_\rho,+_w,\cdot)^\Delta
\end{equation}
be the symmetrization of $(\bar R_\rho,+_w,\cdot)$.
\end{defn}

It is a ring by construction (\cf \S \ref{symmsect}).
We need to ensure that $1\neq 0$ in $W(R,\rho)$. By \eqref{69n} it is enough to show that an equation of the form $1+_wx = x$ where $x\in \bar R_\rho$ gives a contradiction. We have $x\le \rho^n$ for some $n$.
For $x\le \rho^n$ one has $x^{1-\alpha}\ge x\rho^{-n\alpha}$ for $\alpha\in[0,1]$, since $\rho^{n\alpha}\ge x^\alpha~\Rightarrow~\rho^{n\alpha} x^{1-\alpha}\ge x$. One has by construction
\[
1+_wx \ge \rho^{S(\alpha)}x^{1-\alpha}\qqq \alpha\in [0,1]
\]
 and hence using $x^{1-\alpha}\ge x\rho^{-n\alpha}$ one gets
\[
1+_wx \ge \rho^{S(\alpha)-n\alpha}x\qqq \alpha\in [0,1].
\]
Take $\alpha = \frac{e^{-n}}{1+e^{-n}}$ one has $S(\alpha)-n\alpha = \log(e^{-n}+1)>0$ and one gets
\begin{equation}\label{71n}
x\le \rho^n~\Rightarrow~ 1+_wx \ge \rho^\beta x,~\beta = \log(e^{-n}+1)
\end{equation}
Thus the equality $1+_wx = x$ implies that $x \geq \rho^\beta x$ for some $\beta>0$. Iterating one gets $x \geq \rho^{m\beta} x$ for all $m\in\N$. But $x\in [\rho^{-n},\rho^n]$ for some $n$ and one gets
$$
\rho^n\geq x\geq \rho^{m\beta} x\geq \rho^{m\beta} \rho^{-n}\qqq m\in \N
$$
which gives a contradiction provided $\rho\neq 1$.

We  define $r(s)\in W(R,\rho)$ for $s\in\R$ by
\begin{equation}\label{73n}
r(0) = 0,~r(s) = (\rho^{\log(s)},0),~s\ge 0,~ r(s) = (0,\rho^{\log|s|}),~s<0
\end{equation}

\begin{thm}\label{thmalg} Assume $\rho\neq 1$.
Then $W(R,\rho)$ is an algebra over $\R$.
\end{thm}
\proof
It is enough to show that $r$ is an injective homomorphism of the field $\R$ of real numbers in the ring $W(R,\rho)$. One has: $r(s_1s_2)=r(s_1)r(s_2)$ using $\log|s_1s_2|=\log |s_1| + \log |s_2|$ and the rule of signs. For $s_j>0$, one has: $r(s_1)+r(s_2)=r(s_1+s_2)$ by Corollary~\ref{semiringexist1}. This extends to arbitrary signs by adding $r(t)$ to both terms. Moreover we have shown above that $r(1)=1\neq 0$ so that we get a homomorphism of $\R$ in $W(\R,\rho)$. It is necessarily injective since $\R$ is a field.
\endproof

\begin{example}\label{ex}{\rm We take  Example~\ref{mainex}, $R=\tilde C(X,(0,\infty))$ for $X$ a compact space. We write
\begin{equation}\label{75n}
\rho(x) = e^{T(x)}\qquad T(x)\ge 0
\end{equation}
One then gets, using \eqref{59n}, the following formula for the addition $+_w$ of  two functions
\begin{equation}\label{76n}
(f+_w g)(x)=(f(x)^{\beta(x)}+g(x)^{\beta(x)})^{T(x)},\qquad\beta(x)=1/T(x)
\end{equation}
}
\end{example}
\begin{prop} Let $\Omega = \{x\in X, \rho(x)>1\}$ (\ie  the complement of $\{x,T(x)=0\}$). Then  the map $f\to f^\beta$ extends to an isomorphism
\begin{equation}\label{77n}
W(R,\rho)\simeq C_b(\Omega,\R)
\end{equation}
\end{prop}
\proof For $\rho^{-n}\le f\le \rho^n$ one has $e^{-n}\le f^\beta\le e^n$. Conversely,  any $h\in C(\Omega,(e^{-n},e^n))$ extends by $f = h^{T(x)}$ to $X$, with $f(x)=1$ for $x\notin\Omega$. This is enough  to show  that  the symmetrization  process gives $C_b(\Omega,\R)$ starting with functions $h(x)$, $x\in\Omega$ such that $h(x)\in [e^{-n},e^n]$ for some $n$.
\endproof

\subsection{The Banach algebra $\overline W(R,\rho)$}\label{Banach}

 The above Example \ref{ex} suggests to define a semi-norm on $W(R,\rho)$ starting from
\begin{equation}\label{78n}
||f||=\text{Inf}\{\lambda\in \R_+\mid f\le \rho^{\log\lambda}\}\qqq f\in \bar R_\rho
\end{equation}
which is finite for elements of $\bar R_\rho$. One has
\begin{equation}\label{prodin}
    ||fg||\le ||f||||g||
\end{equation}
 because $a\le c$, $b\le d~\Rightarrow~ab\le cd$ in $\bar R_\rho$. Note also that, since $f_1^\alpha\le f_2^\alpha$ if $f_1\le f_2$ one gets
\begin{equation}\label{79n}
f_1\le f_2~\Rightarrow~f_1 +_w g\le f_2 +_w g
\end{equation}
Since $\rho^{\log\lambda_1}+_w\rho^{\log\lambda_2}=\rho^{\log(\lambda_1+\lambda_2)}$ one thus gets
\begin{equation}\label{80n}
||f+_w g||\le||f||+||g||
\end{equation}
We now extend $||.||$ to the symmetrization $W(R,\rho)$. Starting with a semi-ring $A$ and a semi-norm $f\to ||f||$ such that \eqref{prodin} and
 \eqref{80n}   hold, one can endow the semi-ring  $A[\Z/2\Z]$ (\cf \S \ref{symmsect}) with the semi-norm
\begin{equation}\label{82n}
||(a,b)||_1=||a||+||b||
\end{equation}
One has the compatibility with the product
\begin{equation}\label{prod1}
||(a,b)(c,d)||_1\leq ||(a,b)||_1||(c,d)||_1
\end{equation}
since the left hand side is
 $||ac+_w bd||+||ad+_w bc||\le ||a||||c||+||b||||d||+||a||||d||+||b||||c||=||(a,b)||_1||(c,d)||_1$.
 We now take the quotient by the ideal $J = \{(a,a), a\in A\}$ as in \S \ref{symmsect}. Thus the quotient semi-norm  is
\begin{equation}\label{83n}
||(a,b)||_1 =\text{Inf}\{||(x,y)||_1,\exists u\in A, a+_w y+_w u = x +_w b+_w  u\}
\end{equation}
It still satisfies \eqref{80n} and \eqref{prod1}.

 We apply the above discussion to $A=W(R,\rho)$. To show that, after completion for the semi-norm
  \eqref{83n}, we get a Banach  algebra $\overline W(R,\rho)$ over $\R$ we still need to show that  the quotient  norm  does not vanish.

\begin{lem} \label{zeroone} Assume $\rho\neq 1$. Let $(a,b)$ be equivalent  to $(1,0)$ modulo $J$, then $||a||\ge 1$.
\end{lem}
\proof Since $(a,b)$ is equivalent  to $(1,0)$  there exists  $c\in \bar R_\rho$ such that
\begin{equation}\label{85n}
a+_w c  = 1+_w b +_w  c
\end{equation}
thus by \eqref{79n}  one has
\begin{equation}\label{86n}
1+_w c\le a+_w c
\end{equation}
Let us assume  that $||a||< 1$. Then $a\le \rho^{-s}$ for some $s>0$, and one  can then find  $t>0$ such that
\[
e^{-s}+e^{-t}=1
\]
which implies, by \eqref{addrho},  that $\rho^{-s}+_w\rho^{-t}=1$. We have, by \eqref{79n},
\[
1+_w c \le a+_w c\le \rho^{-s}+_w c
\]
and thus
\begin{equation}\label{87n}
\rho^{-t}+_w \rho^{-s} +_w c\le \rho^{-s} +_w c
\end{equation}
With $x =\rho^{-s}+_w c$ this gives  $\rho^{-t}+_w x\le x$ and $1+_w\rho^t x\le \rho^t x$. Since
 $1+_w\rho^t x\geq \rho^t x$ by \eqref{79n}, one gets $1+_w\rho^t x= \rho^t x$ which shows that $1=0$ in $W(R,\rho)$ and contradicts
 Theorem  \ref{thmalg}.
\endproof
It follows that the Banach algebra obtained  as the completion of $W(R,\rho)$ is  non trivial
since the norm of $(1,0)$ is equal to $1$.
Note that  with the notation \eqref{73n} one has
\begin{equation}\label{88n}
||r(s)||=|s|\qquad\forall s\in\R
\end{equation}
Indeed one has $||r(s)||\le |s|$ by construction  using \eqref{73n} and \eqref{78n} and by Lemma \ref{zeroone} one has $||r(1)||=1$ so that  by \eqref{prod1} one has
\[
||r(s)r(1/s)||\le ||r(s)||||r(1/s)||
\]
and $||r(s)||\ge |s|$.
\begin{thm} \label{thmbanach} Assume $\rho\neq 1$. The completion $\overline W(R,\rho)$ is a unital Banach algebra over $\R$.

For any character $\chi$ of the complexification $\overline W(R,\rho)_\C=\overline W(R,\rho)\otimes_\R\C$ one has
\begin{equation}\label{89n}
\chi(\rho)=e
\end{equation}
\end{thm}
\proof  By construction  $\overline W(R,\rho)$ is a real  unital Banach  algebra. The norm on the complexification  is defined by
\begin{equation}\label{90n}
||a+ ib|| = ||a||+||b||
\end{equation}
For a character  $\chi$ of $\overline W(R,\rho)_\C$, the restriction  to $\C 1$ is the identity and thus one has  $\chi(r(s))=s$ which gives \eqref{89n} taking $s=e=2,71828...$
\endproof

It follows from Gelfand's theory that the characters of the complex Banach algebra $\overline W(R,\rho)_\C$ form a non-empty compact space
\begin{equation}\label{spec}
    X=\Spec(\overline W(R,\rho)_\C)\neq \emptyset
\end{equation}
We shall not pursue further the study of this natural notion of spectrum associated to the original pair $(R,\rho)$.

\begin{rem}{\rm Example~\ref{ex} shows that  the norm $||.||_1$ in $W(R,\rho)$ is not the spectral radius norm. In that example it gives $||f||_1=||f_+||+||f_-||$ where $f = f_++f_-$ is the decomposition  of $f$ in positive  and negative parts. The ordinary $C(X,\R)$ norm is given  by $\text{Sup}(||f_+||,||f_-||)$. From \eqref{78n} one gets  in general
\begin{equation}\label{84n}
||f^n||=||f||^n\qquad\forall f\in \bar R_\rho
\end{equation}
However  this holds only before  passing to  the symmetrization.
}\end{rem}
\section{Towards $\run$}\label{sectrun}

Our goal in this section is to show how to apply the above analogue of the Witt construction to the semi-field $\rmax$ of idempotent analysis \cite{Maslov} \cite{Litvinov} and show that it gives in that case the inverse operation of the ``dequantization". It will allow us to take a first speculative step towards the construction of the sought for ``unramified" extension $\run$ of $\R$.

In the case of the Witt construction, the functoriality allows one to apply the functor $\W_p$ to an algebraic closure $\bar{\F}_p$ of $\F_p$ which yields the following diagram
\begin{equation}\label{1o}
\renewcommand{\arraystretch}{1.3}
\begin{array}[c]{ccccc}
\bar{\F}_p&\stackrel{\W_p}{\rightarrow}& \W_p(\bar\F_p)\\
\rotatebox{90}{$\subset$}&&\rotatebox{90}{$\subset$}&&\\
\F_p&\stackrel{\W_p}{\rightarrow}&\Z_p=\W_p(\F_p)
\end{array}
\end{equation}
In our case,  the analogue of the extension $\F_p\subset \bar\F_p$ is the extension of semi-rings
\begin{equation}\label{2o}
\B\subset \rmax
\end{equation}
and the one-parameter group of automorphisms $\theta_\lambda\in\text{Aut}(\rmax)$, $\theta_\lambda(x)=x^\lambda$, plays the role of the Frobenius. But since our construction of $W(\rmax,\rho)$ depends upon the choice of $\rho$, one  first needs to eliminate the choice of $\rho$ by considering simultaneously all possible choices.

\subsection{The $w(\alpha,T)$}

To eliminate the dependence on $\rho$ it is natural  to allow all values of $\rho$, \ie to introduce a parameter $T\ge 0$,
\begin{equation}\label{3o}
\rho=e^T\in \rmax,\qquad\rho\ge 1
\end{equation}
With this notation, $w(\alpha)$ depends on $T$ as it does in the Witt case, one has explicitly
\begin{equation}\label{5o}
w(\alpha,T) = \alpha^{-T\alpha}(1-\alpha)^{-T(1-\alpha)}
\end{equation}
We view the $w(\alpha,T)$ as the analogues of the $w_p(\alpha,T)$ of \eqref{wp} of the Witt case.
The presence in $w(\alpha,T)$ of the parameter $T\geq 0$ means that even if one adds terms which are independent of $T$ the result will depend on $T$.  Thus one works with functions  $f(T)\in\rmax$ with the usual pointwise product and the new addition
\begin{equation}\label{4oo}
(f_1 +_w f_2)(T)=\sum_{\alpha\in\bar I}w(\alpha,T)f_1(T)^\alpha f_2(T)^{1-\alpha}
\end{equation}

\begin{lem}\label{expliadd}
The addition \eqref{4oo} is given by
\begin{equation}\label{4o}
(f_1+_w f_2)(T)=(f_1(T)^{1/T}+f_2(T)^{1/T})^T
\end{equation}
for $T>0$ and by
\begin{equation}\label{4ot0}
(f_1+_w f_2)(0)=\sup(f_1(0),f_2(0))
\end{equation}
\end{lem}
\proof This follows from Lemma \ref{tlem} for $T>0$ and from the equality $w(\alpha,0)=1$ for $T=0$. One then uses Lemma \ref{ineqpow} to conclude in the case $T=0$.
\endproof

\begin{cor}\label{addcons}
The sum of $n$ terms $x_j$ independent of $T$  is given by
\begin{equation}\label{sumofnbis}
    x_1+_w\cdots +_w x_n=\left(\sum x_j^{1/T}\right)^T
\end{equation}
\end{cor}
In particular one can compute the sum of $n$ terms all equal to $1$ which will necessarily be fixed under any automorphism of the obtained structure. One gets
\begin{equation}\label{1ntimes}
    1+_w1+_w\cdots +_w1=n^T
\end{equation}
We expect more generally that the functions of the form
$$
f(T)=x^T\qqq T\geq 0
$$
will be fixed by the lift of the $\theta_\lambda\in \Aut(\rmax)$. We now review the analogy with the Witt case.

\subsection{\te lift}

The constant functions $T\mapsto x$ are the analogue of the Teichm\"uller representatives
\begin{equation}\label{11o}
\tau(x)(T) = x\qquad\forall T
\end{equation}
One has
\begin{equation}\label{12o}
\tau(x)+\tau(y) =  \sum_{\alpha\in \bar I}w(\alpha,T)\,x^\alpha y^{1-\alpha}
\end{equation}
where the sum in the right hand side is computed  in $\rmax$. We view this formula  as the analogue  of the formula \eqref{24m} of the usual Witt case.

\subsection{Residue morphism}\label{resmorhismsect}

The evaluation at $T=0$  is by construction a morphism
\begin{equation}\label{9o}
\epsilon: f\mapsto f(0)\in \rmax
\end{equation}
We view this morphism as the analogue of the canonical map which exists for any strict $p$-ring
\begin{equation}\label{10o}
\epsilon_p:\W_p(K) \to K=\W_p(K) /p\W_p(K)
\end{equation}

\subsection{Lift of the automorphisms $\theta_\lambda\in\Aut(\rmax)$}

One has a natural one parameter group of automorphisms $\alpha_\lambda$ of our structure, which corresponds to the $\theta_\lambda\in\Aut(\rmax)$. It is given  by
\begin{prop}\label{autom} The following defines a one parameter group of automorphisms
\begin{equation}\label{6o}
\alpha_\lambda(f)(T) = f(T/\lambda)^\lambda\qquad\forall \lambda\in\R_+^\times
\end{equation}
One has
\begin{equation}\label{rescom}
   \epsilon\circ\alpha_\lambda = \theta_\lambda\circ\epsilon \qqq \lambda
\end{equation}
and
\begin{equation}\label{teichcom}
    \alpha_\lambda\circ \tau=\tau\circ \theta_\lambda \qqq \lambda.
\end{equation}
\end{prop}
\proof One has for $T>0$ using Lemma \ref{expliadd},
\begin{align*}
\alpha_\lambda(f_1+_w f_2)(T)&=\left( (f_1+_w f_2)(T/\lambda)\right)^\lambda = \left(f_1(T/\lambda)^{\lambda/T}+f_2(T/\lambda)^{\lambda/T}\right)^{(T/\lambda)\times \lambda}\\&=
(\alpha_\lambda(f_1)(T)^{1/T} + \alpha_\lambda(f_2)(T)^{1/T})^T =(\alpha_\lambda(f_1)+_w \alpha_\lambda(f_2))(T)
\end{align*}
thus $\alpha_\lambda$ is additive. It is also multiplicative and defines an automorphism.
The equality \eqref{rescom} follows from \eqref{6o} by evaluation at $T=0$.
The equality \eqref{teichcom} also  follows from \eqref{6o}. \endproof

\subsection{Fixed points}

We now determine the fixed points of the  $\alpha_\lambda$ and show, as expected from \eqref{1ntimes} that they form the semi-field $\R_+$.

\begin{prop} \label{fixedpts}
The fixed points of $\alpha_\lambda$ are of the form
\begin{equation}\label{7o}
f(T) = a^T
\end{equation}
and they form  the semi-field $\R_+$ which is the positive part of the field $\R$ of real numbers,
endowed with the ordinary addition.
\end{prop}

\proof Assume that $\alpha_\lambda(f)=f$ for all $\lambda>0$. Then one has, using \eqref{6o} for $\lambda=T$, $f(T)=f(1)^T$ which gives \eqref{7o}. Moreover by Lemma \ref{expliadd} the addition corresponds to the  semi-field $\R_+$ which is the positive part of the field $\R$ of real numbers.
\endproof

\subsection{Characters and representation by functions}\label{charep}

For each $T>0$ the algebraic operations on the value $f(T)$ are the same as in the semi-field $\R_+$ using the evaluation $f(T)^{1/T}$. Thus there is a uniquely associated character $\chi_T$ which is such that
\begin{equation}\label{charach}
    \chi_T(f)=f(T)^{1/T}
\end{equation}
and we now use the characters $\chi_T$ to represent the elements of the extension $\run$ as functions of $T$ with the ordinary  operations of pointwise sum and product.
\begin{prop}\label{funrepr}
The following map $\chi$ is a homomorphism of semi-rings to the algebra of functions from $(0,\infty) $ to $\R_+$ with pointwise sum and product,
\begin{equation}\label{chimap}
    \chi(f)(T)=f(T)^{1/T}\qqq T>0.
\end{equation}
One has
\begin{equation}\label{teichrep}
    \chi(\tau(x))(T)=x^{1/T}\qqq T>0
\end{equation}
and
\begin{equation}\label{autrep}
   \chi(\alpha_\lambda(f))(T)=\chi(f)(T/\lambda)\qqq T>0.
\end{equation}
\end{prop}
\proof These properties are straightforward consequences of \eqref{charach}.\endproof

In this representation $\chi$, the residue morphism of \S \ref{resmorhismsect} is given, under suitable continuity assumptions by
\begin{equation}\label{reschi}
    \epsilon(f)=\lim_{T\to 0}\chi(f)(T)^T\,.
\end{equation}
In this representation the algebraic operations are very simple and this suggests to represent elements of $\run$ as  functions $\chi(f)(T)$. Among them one should have the fixed points \eqref{7o} which give $\chi(f)=a$ and the \te lifts which give \eqref{teichrep}. We parameterize the latter in the form
\begin{equation}\label{exps}
    e_\xi(T)=e^{-\xi/T}\qqq T>0.
\end{equation}
After symmetrization and passing to the field of quotients, the fixed points  \eqref{7o} and the \te lifts \eqref{11o} generate the field of fractions of the form (in the $\chi$ representation)
\begin{equation}\label{fracchi}
 \chi(f)(T)=  \left(\sum a_j e^{-\xi_j/T}\right)\text{\huge/}\left(\sum b_j e^{-\eta_j/T}\right)
\end{equation}
where the coefficients $a_j,b_j$ are real numbers and  the exponents $\xi_j,\eta_j\in \R$. While such expressions give a first hint towards $\run$ one should not be satisfied yet since, as explained in \S \ref{deforsect} below, natural examples coming from quantum physics use expressions of the same type but involving more elaborate sums. In all these examples, including those coming from the functional integral, it turns out that not only $\lim_{T\to 0}\chi(f)(T)^T$ exists as in \eqref{reschi} but in fact the  function $f(T)=\chi(f)(T)^T$ admits an asymptotic expansion for $T\to 0$ of the form
\begin{equation}\label{asexp}
    f(T)=\chi(f)(T)^T\sim \sum a_nT^n\,.
\end{equation}
For functions of the form \eqref{fracchi}, this expansion only uses the terms with the lowest values of $\xi_j$ and $\eta_j$ and is thus only a crude information on the element $f$. Moreover as soon as one uses integrals instead of finite sums in \eqref{fracchi} one obtains general asymptotic expansions \eqref{asexp}. In fact, given a convergent series
$
g(T)=\sum_0^\infty b_nT^n
$
one has (\cf \eg \cite{Ramis})
\begin{equation}\label{borel}
g(T)=b_0+\int_0^\infty e^{-\xi/T}\phi(\xi)d\xi
\end{equation}
where $\phi(\xi)$ is the Borel transform
$$
\phi(\xi)=\sum_0^\infty b_{n+1}\frac{\xi^n}{n!}
$$
Thus $g(T)$ is of the form \eqref{fracchi} using integrals. Moreover, with $b_0=1$, the asymptotic expansion of $g(T)^T$ for $T\to 0$ is of the form
\begin{equation}\label{asexpbis}
    g(T)^T\sim \sum a_nT^n\,,
\end{equation}
where the coefficients $a_n$ are given by $a_0=1$ and
$$
\begin{array}{cc}
 a_1 =& 0 \\
 a_2 =& b_1 \\
 a_3 =& -\frac{b_1^2}{2}+b_2 \\
 a_4 =& \frac{b_1^2}{2}+\frac{b_1^3}{3}-b_1 b_2+b_3 \\
 a_5 =& -\frac{b_1^3}{2}-\frac{b_1^4}{4}+b_1 b_2+b_1^2 b_2-\frac{b_2^2}{2}-b_1 b_3+b_4 \\
 a_6 =& \frac{b_1^3}{6}+\frac{11 b_1^4}{24}+\frac{b_1^5}{5}-\frac{3}{2} b_1^2 b_2-b_1^3 b_2+\frac{b_2^2}{2}+b_1 b_2^2+b_1 b_3+b_1^2 b_3-b_2 b_3-b_1 b_4+b_5
\end{array}
$$
which shows how to determine the $b_n$ once the coefficients $a_j$, $j\leq n+1$ are given. Using the freedom to multiply $g(T)$ by $ae^{-\xi_0/T}$ for $a,\xi_0\in \R$, $a>0$, one gets in this way any convergent series $\sum a_nT^n$ with $a_0>0$ as the asymptotic expansion of $g(T)^T$ for $g(T)$ of the form
$$
g(T)=ae^{-\xi_0/T}+\int_{\xi_0}^{\xi_1} e^{-\xi/T}\phi(\xi)d\xi
$$
One checks that the only restriction on the $a_n$ is that $a_0>0$ so that $a_0=e^{-\xi_0}$ for some $\xi_0\in \R$.
This suggests more generally to use the theory of divergent series (\cf \cite{Ramis}) in the construction of $\run$. The simple point is that the action of $\R_+^*$ given by the $\alpha_\lambda$ of \eqref{6o}, gives a grading which admits the $f_n$, $\chi(f_n)(T)=T^n$ as eigenvectors. This is only formal since the $f_n(T)=T^{nT}$ do not have an asymptotic expansion \eqref{asexpbis} but the $\chi(f_n)$ are integrals of the above form since
$$
n!\,T^{n+1}=\int_0^\infty e^{-\xi/T}\xi^n d\xi\,.
$$

\subsection{Deformation parameter $T\sim\hbar$}\label{deforsect}
 In idempotent analysis (\cite{Maslov}, \cite{Litvinov}), the  process that allows one to view $\rmax$ as a result of a deformation of the usual algebraic structure on real numbers is known as ``dequantization'' and can be described as a semi-classical limit in the following way.  First of all note that, for real numbers $w_j$,   one has
\[
 h\ln(e^{w_1/h}+e^{w_2/h}) \to \max\{w_1,w_2\}\quad\text{as}~h\to 0.
\]
Thus,  the natural map
\[
D_h: \R_+ \to \R_+\qquad D_h(u) = u^h
\]
satisfies the equation
$$
\lim_{h\to 0}D_h(D_h^{-1}(u_1)+D_h^{-1}(u_2))=\max \{u_1,u_2\}.
$$
In the limit $h\to 0$, the usual algebraic rules on $\R_+$ deform to become those of $\rmax$.
In our context this corresponds to the residue morphism of \S \ref{resmorhismsect} expressed in the $\chi$-representation. More specifically, with $\chi(f)$ given by \eqref{fracchi} one has
\begin{equation}\label{limchi}
    \epsilon(f)=\lim_{T\to 0}\chi(f)(T)^T=e^{-\inf_j(\xi_j)+\inf_j(\eta_j)}
\end{equation}
when the coefficients $a_j,b_j$ are positive real numbers so that the exponentiation makes sense. These elements form a semi-field and $\epsilon$ is a
homomorphism  to the semi-field $\rmax$. 
In fact, it extends to a homomorphism  $\tilde\epsilon$ from the field $K$ of rational fractions of the form \eqref{fracchi} to the hyperfield $\cT\R$ of tropical reals defined by O. Viro (\cf \cite{viro}) by the following hyperaddition of real numbers:
\begin{equation}\label{taur}
    a \smile b=\left\{
                 \begin{array}{ll}
                   a, & \hbox{if $|a|>|b|$ or $a=b$;} \\
                   b, & \hbox{if $|a|<|b|$ or $a=b$;} \\
                   $[-a,a]$, & \hbox{if $b=-a$.}
                 \end{array}
               \right.
\end{equation} The extension is given by $\tilde\epsilon(0)=0$ and for a reduced fraction
\begin{equation}\label{fracchi1}
 \chi(f)(T)=  \left(\sum a_j e^{-\xi_j/T}\right)\text{\huge/}\left(\sum b_j e^{-\eta_j/T}\right)
\end{equation}
\begin{equation}\label{troprun1}
    \tilde\epsilon(f)={\rm sign}(\frac{a_{j_0}}{b_{k_0}})e^{-\xi_{j_0}+\eta_{k_0}}\,, \ \xi_{j_0}=\inf_j(\xi_j), \ \eta_{k_0}=\inf_k(\eta_k).
\end{equation}

One has
 $
    \chi(f)(T)\sim a\, \tilde\epsilon(f)^{1/T}\,, \ \text{for}\, \, T\to 0,
 $
 for some $a> 0$ (with the notation $x^\lambda={\rm sign}(x)|x|^\lambda$ for $x\in \R$ and $\lambda>0$). Thus one obtains the hyperfield  $\cT\R$ of tropical reals as the quotient
 $K/G$ of the field $K$ by a subgroup of its multiplicative group which, as already observed by M. Krasner, is the natural construction of many hyperfields (\cite{wagner}, \cite{japanese}).
 
The formalism of idempotent analysis, motivated by quantum physics,  suggests that the parameter $T$ should be related to the Planck constant $\hbar$. Moreover in order to use $\run$ in the context of quantum physics,  one should relax the requirement that the sums involved in \eqref{fracchi} only involve finitely many terms. The key example  is given by the functional integral in the Euclidean formulation of Quantum Field Theory. Indeed the generating function of Euclidean Green's
functions is given by (\cf \eg \cite{CMbook})
\begin{equation}\label{ZJEbis}
 Z(J_E) = \cN \int \exp\left(-\frac{S(\phi_E)- \langle J_E, \phi_E
 \rangle}{\hbar}\right)\,
 \cD[\phi_E]
\end{equation}
where $S(\phi_E)$ is the Euclidean action, in terms of the Euclidean classical fields $\phi_E$, the source  $J_E$ is an element of the linear space dual
to that of Euclidean classical fields and the normalization factor $\cN $ is the inverse of
$$
\int \exp\left(-\frac{S(\phi_E)}{\hbar}\right)\,
 \cD[\phi_E]
 $$
Such integrals are typical sums involving $+_w$ where $w$ is the function  given by \eqref{5o}
but since the sums are infinite one needs to extend the entropy from finite partitions of $1$ to infinite partitions. The formula for computing the sum is then deeply related to the  basic formula of thermodynamics using the entropy to express the free energy from a variational principle involving the sum of the entropy and the energy with suitable multipliers. It suggests that it might be worthwhile to reconsider the functional integral from this angle, considering this formula as more basic than the analogy with ordinary integrals.

Note finally that the expansion \eqref{asexp} still holds for the elements such as $Z(J_E)$ of \eqref{ZJEbis} and that this asymptotic expansion is in general much more involved than in the simplest example of \eqref{fracchi} since it is the ``loop expansion" of quantum field theory which is the basis of the concrete computations in quantum physics (\cf \eg \cite{CMbook}). In conclusion the above development suggests that the extension $\run$ of $\R$ is the proper receptacle for the ``values" of many $\hbar$-dependent physical quantities arising in quantum field theory. Together with the previous understanding of renormalization from the Riemann-Hilbert correspondence (\cf \cite{CK}, \cite{cknew}, \cite{cmln}, \cite{CMbook}) this should be an important piece of the puzzle provided by the quantum.

\end{document}